\documentclass[]{article}
\usepackage[utf8]{inputenc} 
\usepackage[T1]{fontenc}    
\usepackage{hyperref}       
\usepackage{url}            
\usepackage{microtype}
\usepackage{graphicx}
\usepackage{array}
\usepackage{color} 
\usepackage{tabularx}
\usepackage{amsmath,mathdots,amsthm}
\usepackage{amssymb}
\usepackage{amsfonts}
\usepackage{txfonts}
\usepackage{arxiv}
\usepackage{xcolor}
\usepackage{bm}
\usepackage{soul}
\usepackage{float}
\usepackage{cite}

\hypersetup{
	colorlinks=true,                          
	linkcolor=blue, % equation, section link color
	citecolor=blue, % bib color
	urlcolor=black  % url color if any
}

\newtheorem*{rmrk*}{Remark}

\newtheorem*{definition*}{Definition}

\newtheorem{thm}{Theorem}
\newtheorem*{thm*}{Theorem}

\newcommand{\nrm}[1]{ \left\vert\left\vert #1 \right\vert\right\vert }
\newcommand{\gradphi}{\nabla \phi}
\newcommand{\dphi}{\dot{\phi}}
\renewcommand{\u}{\mathbf{u}}
\newcommand{\m}{\mathbf{m}}
\newcommand{\estar}{\varepsilon^\star}
\renewcommand{\div}[1]{\mathrm{div}\left( #1 \right)}
\newcommand{\pd}[2]{\frac{\partial #1}{\partial #2}}
\newcommand{\prn}[1]{\left( #1 \right)}

\renewcommand{\j}{\mathbf{q}}

\newcommand{\grad}[1]{\nabla\left( #1 \right)}
\newcommand{\dx}{\boldsymbol{\zeta}}
\newcommand{\ppsi}{\boldsymbol{\psi}}
\newcommand{\x}{\mathbf{x}}
\newcommand{\tr}[1]{\mathrm{Tr}\left( #1 \right)}
\newcommand{\blue}[1]{\textcolor[rgb]{0.00,0.40,0.80}{#1}}
\newcommand{\red}[1]{\textcolor[rgb]{0.95,0.20,0.20}{#1}}
\newcommand{\pdt}[1]{\frac{\partial #1}{\partial t}}
\newcommand{\RH}[1]{\left[  #1 \right]}

\newcommand{\dblint}[1]{\int_{t_0}^{t_1} \!\!\! \int_{\Omega_t} \   \ #1\  \ d\Omega \ dt}
\newcommand{\calE}{\mathcal{E}}
\newcommand{\calR}{\mathcal{R}}
\newcommand{\calM}{\mathcal{M}}
\newcommand{\calD}{\mathcal{D}}
\newcommand{\calG}{\mathcal{G}}

\renewcommand{\j}{\mathbf{j}} % \gradphi

\title{An Eulerian  hyperbolic  model for heat transfer derived via Hamilton's principle: analytical and numerical study}

\author{
	 \href{https://orcid.org/0000-0002-7150-3313}{\includegraphics[scale=0.06]{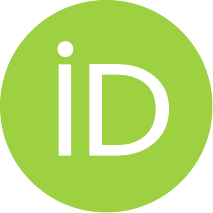}\hspace{1mm}Firas Dhaouadi}\thanks{corresponding author.\newline Email addresses: \texttt{firas.dhaouadi@unitn.it} (F. Dhaouadi), \texttt{sergey.gavrilyuk@univ-amu.fr} (S. Gavrilyuk)} $^{\;a}$, \hspace{3mm}
	\href{https://orcid.org/0000-0003-4605-8104}{\includegraphics[scale=0.06]{orcid.pdf}\hspace{1mm}Sergey Gavrilyuk$^b$} 
}

\affiliation{$^a$ Department of Civil, Environmental and Mechanical Engineering, University of Trento, Via Mesiano 77, 38123 Trento, Italy. \\ 
$^b$Aix-Marseille University and CNRS UMR 7343 IUSTI, 5 rue Enrico Fermi, 13453 Marseille, France.
}
\date{}

\renewcommand{\headeright}{F. Dhaouadi AND S. Gavrilyuk}
\renewcommand{\undertitle}{}
\renewcommand{\shorttitle}{}

\hypersetup{
pdftitle={An Eulerian  hyperbolic  model for heat transfer derived via Hamilton's principle: analytical and numerical study},
pdfsubject={q-bio.NC, q-bio.QM},
pdfauthor={Firas Dhaouadi, Sergey Gavrilyuk},
pdfkeywords={First keyword, Second keyword, More},
}

\begin{document}
\maketitle

\begin{abstract}
	In this paper, we present a new model for heat transfer in compressible fluid flows. The model is derived from Hamilton's principle of stationary action in Eulerian coordinates, in a setting where the entropy conservation is recovered as an Euler--Lagrange equation. The governing system is shown to be hyperbolic. It is asymptotically consistent with the Euler equations for compressible heat conducting fluids, provided the addition of suitable relaxation terms. A study of the Rankine--Hugoniot conditions and the Clausius--Duhem inequality reveals that contact discontinuities cannot exist while expansion waves and compression fans are possible solutions to the governing equations. Evidence of these properties is provided on a set of numerical test cases.
\end{abstract}

% keywords can be removed

\textit{\red{The final version of this paper has been published in Proceedings of the Royal Society A, Volume 480, Issue 2283, 2024. The published version is available at} \blue{https://royalsocietypublishing.org/doi/abs/10.1098/rspa.2023.0440} }

	\section{Introduction}
The description of heat transfer processes in continuum mechanics is generally based on  Fourier's law of heat conduction stating that the heat flux  in a medium is opposite  to the temperature gradient.  However, one flaw in this description is the fact that a localized perturbation in the temperature field, generates an instantaneous response, as small as it may be, in all the medium. There have been several attempts to circumvent this shortcoming, maybe the most known of which is the Maxwell-Cattaneo-Vernotte model of  heat conduction  \cite{cattaneo1948,maxwell2003dynamical,vernotte1958paradoxes},  where the heat equation was approximated by a hyperbolic system of evolution equations for the temperature and the heat flux. Such a coupling between the heat processes and  motion of continuum media, where the heat flux is considered as an independent variable was considered in particular in  \cite{cimmelli1991new,muller2013rational,greennaghdi1991,green1993thermoelasticity,peshkov2018continuum}.

In this context, we present in this paper a new  model describing heat transfer in an inviscid fluid flow. The model is derived  via Hamilton's principle of stationary action in Eulerian coordinates resulting in a first order quasilinear system of partial differential equations. 
The formulation we propose finds its foundation in the works of Green \& Naghdi \cite{greennaghdi1991,green1993thermoelasticity} and allows to recover the entropy evolution equation as an Euler-Lagrange equation for the Hamilton action. Under the classical assumption of the total energy convexity, the governing equations are shown to be hyperbolic. The system admits two types of waves which we call \textit{acoustic} and \textit{thermal} waves. The former is analogous to the acoustic waves of the Euler equations while the latter describes  rather the heat transfer. We show for acoustic waves that the corresponding eigenfields are genuinely nonlinear in the sense of Lax \cite{lax1973hyperbolic},  while the eigenfields for thermal waves are neither genuinely linear nor linearly degenerate. In particular, in the last case we show the existence of {\it expansion shock}s, \textit{i.e.}, discontinuous fronts moving through the media for which the pressure behind the shock is lower than the pressure ahead of it, and related {\it compression fans} as possible solutions of the governing equations.  The theory of expansion shocks was developed in the case of a scalar quasilinear conservative equation  with a non-convex flux in \cite{oleinik1959uniqueness} and extended further to systems of hyperbolic equations in  \cite{wendroff1972riemann,wendroff1972riemann2,liu1981admissible,lefloch2002hyperbolic}.  Physical applications of expansion shocks are described and discussed for example in \cite{zel1946, dettleff1979experimental, borisov1983rarefaction,kutateladze1987rarefaction, menikoff1989riemann,kluwick2018shock}.

The rest of this paper is organized as follows. In section \ref{sec:Euler}, we recall the classical governing equations for compressible fluid flows accounting for heat conduction. 
Section \ref{sec:model} provides the theoretical foundation for our proposed model. It contains details on the  system Lagrangian, the governing equations, and how dissipation is added by means of the Rayleigh dissipation function, in a way that is asymptotically compatible with the Fourier law. We also provide in this section a criterion for the total energy convexity and present our equations as a symmetric $t$-hyperbolic system in the sense of Friedrichs.  In section \ref{sec:Hyperbolic}, we provide an extensive study of the model hyperbolicity by analyzing its eigenstructure and investigating the nature of the corresponding eigenfields. We also draw a comparative analysis of the dispersion relations of our model with the Euler system with heat conduction where we show that both models are compatible in the low frequencies up to a certain threshold. In section \ref{sec:shocks}, we prove the non-existence of contact discontinuities for our system and we show that there are two distinct branches of the Hugoniot curve, of acoustic and thermal natures. In particular, we use the Clausius-Duhem inequality  as a selection criterion for admissible shocks and we show that expansion shocks solutions exist in certain conditions. Section \ref{sec:Numerics} lays out all the necessary details on the implemented numerical scheme used to numerically solve the governing system of equations. Lastly, we present in section \ref{sec:results} a sample of numerical tests for specific Riemann problems, confirming the theoretical findings on the nature of solutions. The paper is supported by three appendixes which contain supplementary details. Appendix \ref{app:calculus} complements section \ref{sec:model} with all the necessary calculus for the model derivation through Hamilton's principle. Appendix \ref{app:energy_conservation} gives the proof of energy conservation in the general case, even when curl-cleaning is applied. Lastly, Appendix \ref{app:sym} contains the proof of symmetrizability of the governing equations.   

\subsection*{Notations}
Tensor notation will be used throughout the paper. Thus, it seems of importance to briefly recall the definitions and conventions we use here for the sake of clarity. For any $\lambda\in\mathbb{R}$,   $\u,\mathbf{v}\in\mathbb{R}^d$ and $\mathbf{A}\in\mathbb{R}^{d\times d}$, we convene what follows. The transpose of $\u$ shall be denoted as $\u^T$ and the trace of $\mathbf{A}$ will be denoted  as $\mathrm{Tr}(\mathbf{A})$. The identity matrix will be denoted as $\mathbf{I}$. We define the differential operators as follows 
\begin{align*}
	\div{\u} = \pd{u_1}{x_1}+\cdots+\pd{u_d}{x_d}, \quad \pd{\lambda}{\x}=\left(\  \pd{\lambda}{x_1},\ \cdots,\ \pd{\lambda}{x_d} \ \right), \quad \nabla\lambda = \left(\pd{\lambda}{\x}\right)^T =\left(\  \pd{\lambda}{x_1},\ \cdots,\ \pd{\lambda}{x_d} \ \right)^T,
\end{align*}
\begin{align*}
	\pd{\u}{\x}  =\left(
	\begin{array}{ccc}
		\displaystyle  \pd{u_1}{x_1}& \cdots & \displaystyle \pd{u_1}{x_d} \\ 
		\vdots & \ddots & \vdots \\
		\displaystyle  \pd{u_d}{x_1}& \cdots & \displaystyle \pd{u_d}{x_d} 
	\end{array}
	\right),
	\qquad 
	\nabla\u =   \left(\pd{\u}{\x}\right)^T  =\left(
	\begin{array}{ccc}
		\displaystyle  \pd{u_1}{x_1}& \cdots & \displaystyle \pd{u_d}{x_1} \\ 
		\vdots & \ddots & \vdots \\
		\displaystyle  \pd{u_1}{x_d}& \cdots & \displaystyle \pd{u_d}{x_d} 
	\end{array}
	\right).
\end{align*}
The divergence of the second order tensor $\mathbf{A}$ is defined such that for every constant vector field $\mathbf{a}$, we have 
\begin{equation*}
	\div{\mathbf{A}}\cdot\mathbf{a} = \div{\mathbf{A}\mathbf{a}}, 
\end{equation*}
\textit{i.e}, the divergence of the second order tensor $\mathbf{A}$ is the vector whose components are the divergence of each column of $\mathbf{A}$.   
In particular, we recall the following formulas
\begin{align*}
	\mathrm{div}(\mathbf{A}\u)=\mathrm{div}(\mathbf{A})\cdot\u+\mathrm{Tr}\left(\mathbf{A}\frac{\partial\u}{\partial \mathbf{x}}\right), \qquad 
	\mathrm{div}(\u\otimes\mathbf{v})= \mathrm{div}(\u)\ \mathbf{v} + \frac{\partial \mathbf{v}}{\partial \mathbf{x}}\  \u, \qquad
	\mathrm{Tr}(\u\otimes \mathbf{v}\ \mathbf{A}) = \mathbf{v}\cdot \mathbf{A}\u .
\end{align*}
The Hessian matrix of  $\lambda$ with respect to the variables $\mathbf{v}=(v_1,\cdots,v_n)^T$ shall be denoted as $\nabla^2_\mathbf{v}\lambda$ and is the symmetric matrix whose  expression is given by 
\begin{equation*}
	\nabla^2_\mathbf{v}\lambda = \left(
	\begin{array}{ccc}
		\displaystyle  \pd{^2\lambda}{v_1^2}& \cdots & \displaystyle \pd{^2\lambda}{v_1v_n} \\ 
		\vdots & \ddots & \vdots \\
		\displaystyle  \pd{^2\lambda}{v_1v_n}& \cdots & \displaystyle \pd{^2\lambda}{v_n^2} 
	\end{array}
	\right).
\end{equation*}
Lastly, the material derivative, with respect to a velocity field $\u$ is defined as follows
\begin{equation*}
	\dot{\lambda} =  \pdt{\lambda} + \pd{\lambda}{\x} \u         \qquad	\dot{\mathbf{v}} = \pdt{\mathbf{v}} + \pd{\mathbf{v}}{\x}\u
\end{equation*}   

\section{The Euler system for compressible heat conducting fluid flows}
\label{sec:Euler}
Compressible flows of inviscid fluids in the presence of heat transfer can be described by the Euler equations supplemented by the Fourier heat conduction terms
\begin{subequations}
	\begin{alignat}{2}
		& \pdt{\rho} &&+ \div{\rho \u } = 0, \label{eq:EF_rho} \\
		& \pdt{\rho\u}  &&+ \div{ \rho \u \otimes \u +p(\rho,\eta)\mathbf{I}  }= 0, \label{eq:EF_rhou} \\
		&\pdt{E} &&+ \div{E\u+p(\rho,\eta)\u - K \nabla \theta(\rho,\eta)} = 0.
		\label{eq:EF_E}
	\end{alignat}
	\label{eq:Euler-Fourier}
\end{subequations}
Applications of such a model can be found for example in  \cite{le2006compressible}. In this part and for the rest of this paper, $\x \in \mathbb{R}^d$ refers to the vector of space coordinates and $t \in \mathbb{R}_+$ is the time. $K$ is the heat conductivity assumed positive, $\rho=\rho(\x,t)$ is the density field, $\u=\u(\x,t)$ is the velocity field, $\theta(\rho,\eta)$ is the temperature, $p(\rho, \eta )$ is the pressure, and $\eta=\eta(\x,t)$ is the specific entropy. For the sake of lightness, dependence of the variables on space and time will be assumed from now on but omitted in the notation. Lastly, $E$ is the total energy of the system which is expressed as
\begin{equation*}
	E= \frac{1}{2}\rho \nrm{\u}^2 + \rho \varepsilon(\rho,\eta),
\end{equation*}
where $\varepsilon(\rho,\eta)$ is the specific internal energy, which is linked to the temperature $\theta(\rho,\eta)$ and to the pressure $p(\rho,\eta)$ via the Gibbs identity
\begin{equation*}
	\theta d\eta=d\varepsilon+p \  d\prn{\frac{1}{\rho}}.
\end{equation*}
so that we have
\begin{equation*}
	p(\rho,\eta)=\rho^2\pd{\varepsilon}{\rho}, \quad \theta(\rho,\eta)=\pd{\varepsilon}{\eta}.
\end{equation*}
Equations \eqref{eq:Euler-Fourier} describe the conservation of mass, momentum and total energy, respectively. Furthermore, an additional evolution equation for the entropy can be derived as a consequence of this system and which can be cast as the Clausius-Duhem inequality
\begin{equation*}
	\pdt{\rho\eta} + \div{\rho \eta \u - \frac{K}{\theta}\nabla\theta} =  \frac{K}{\theta^2}\nrm{\nabla\theta}^2 \geq0
\end{equation*}
\section{Model description} 
\label{sec:model}
\subsection{Lagrangian of the system}
Derivation of the Euler equations of compressible flows from Hamilton's principle of stationary action can be found in \cite{Berdichevsky2009, godunov2003elements, Gavrilyuk_2011, serre1993principe}. The corresponding Lagrangian can be written as    
\begin{equation*}
	\mathcal{L}_E = \int_{\Omega_t} \ \prn{\frac{1}{2}\rho \nrm{\u}^2  - \rho \varepsilon(\rho,\eta)} \ d\Omega,
\end{equation*}
where $\Omega_t$ is the material domain. 
The momentum equation is then obtained as the Euler-Lagrange equation while the mass and entropy conservation laws are considered as constraints. However, postulating the entropy equation as a constraint is not necessary as it can be recovered as an Euler-Lagrange equation as well (See equation \eqref{eq:EL_1} in Appendix \ref{app:calculus}). Indeed, in a series of papers by Green and Naghdi \cite{greennaghdi1991,green1993thermoelasticity}, an independent auxiliary potential $\phi(\x,t)$ was introduced as a thermal analogue of the kinematic variables. In particular, it is linked to the temperature by the relation
\begin{equation}
	\dphi(\x,t) = -\theta(\x,t)
	\label{eq:scalar_phi}
\end{equation}
A similar idea was also used in Lagrangian coordinates in  \cite{peshkov2018continuum}. In what follows, the governing equations will be derived in Eulerian coordinates. For this, consider the following Lagrangian 
\begin{equation}
	\mathcal{L} = \int_\Omega \ \prn{\frac{1}{2}\rho \nrm{\u}^2  - \frac{1}{2}\alpha(\rho)  \nrm{\gradphi}^2 - \rho \estar(\rho,\dphi)} \ d\Omega,
	\label{eq:Lagrangian}
\end{equation}
where the coefficient $\alpha(\rho)$ is an arbitrary positive function of density.  The quantity $\estar$ is the Helmholtz free energy, which relates to the internal energy by a partial Legendre transform as follows
\begin{equation*}
	\varepsilon(\rho, \eta)= \estar (\rho,\dot \phi)-\eta \dot \phi, \quad {\rm with }\quad \eta=\frac{\partial \estar }{\partial \dot \phi}.
\end{equation*}
In addition, one has
\begin{equation*}
	\frac{\partial \varepsilon }{\partial \rho}=\frac{\partial \estar}{\partial \rho}, \quad \frac{\partial \varepsilon }{\partial \eta}=-\dot \phi.
\end{equation*}
If we recall that $\theta=-\dot \phi $ and note the specific volume by $v=1/\rho$, then the internal energy verifies the Gibbs identity
\begin{equation}
	\theta d\eta=d\varepsilon+pdv.
	\label{eq:Gibbs}
\end{equation}
In order to obtain the governing system of equations, a variational principle is applied to the Lagrangian \eqref{eq:Lagrangian}, submitted to the mass conservation constraint
\begin{equation}
	\pdt{\rho} + \div{\rho \u} = 0.
	\label{eq:mass}
\end{equation} 
By virtue of their definitions, the variations of both the density $\rho$ and the velocity $\u$ can be linked to the kinematic displacement of the continuum while $\gradphi$ and $\dphi$ are derivatives of the \textit{thermal displacement}. This suggests that two Euler-Lagrange equations are to be obtained and which are given hereafter, respectively
	\begin{align*}
		&	\pdt{} \prn{\rho \u} + \div{\rho \u \otimes \u + \prn{\rho^2 \pd{\varepsilon}{\rho}+\frac{1}{2}\prn{\rho\alpha'(\rho)-\alpha(\rho)}\nrm{\nabla\phi}^2}\mathbf{I} + \alpha(\rho) \ \gradphi \otimes \gradphi }= 0, 
		\\
		& 	\pd{}{t}\prn{\rho \eta} + \div{ \rho \eta\u + \alpha(\rho) \gradphi} = 0.
	\end{align*}
Details on the derivation of this system are given in Appendix \ref{app:calculus}. In order to recast these equations into a first-order quasilinear system, we can apply an order reduction by denoting $\j=\gradphi$. Its evolution is obtained by applying the gradient operator on \eqref{eq:scalar_phi} so that 
\begin{equation*}
	\pd{\j}{t} + \grad{\j\cdot\u + \theta} = 0.
\end{equation*}
It is worthy of note in this context that $\j$ is a gradient field, evolved as an independent variable. In virtue of this definition, the curl-free constraint $\nabla\wedge\j=0$ must be respected by admissible solutions. This also allows to restructure the corresponding evolution equation, by adding $\prn{\nabla\wedge\j}\wedge\u$ to the left-hand side. This subtle modification is important, as it restores Galilean invariance to the system. 
Taking into account all the previous considerations and notations, one can obtain the following first-order set of evolution equations  
\begin{subequations}
	\begin{alignat}{2}
		& \frac{\partial \rho}{\partial t} &&+ \div{\rho \u } = 0, \\
		& \frac{\partial \rho\u}{\partial t}  &&+ \div{ \rho \u \otimes \u +\Pi  }= 0, \quad \Pi = P(\rho,\eta, \j) \ \mathbf{I}+\alpha(\rho)\ \j \otimes \j,  \\
		& \frac{\partial \j}{\partial t} &&+ \grad{\j\cdot\u + \theta(\rho,\eta)}  + \prn{ \frac{\partial \j}{\partial \x} - \prn{\frac{\partial \j}{\partial \x}}^T }\u = 0, \\
		&\frac{\partial \rho\eta}{\partial t} &&+ \div{\rho \eta\u+ \alpha(\rho) \j  } = 0.
	\end{alignat}
	\label{eq:hyp_heat_eq}
\end{subequations}
where the total pressure $P(\rho,\eta, \j)$ is given by 
\begin{equation*}
	P(\rho,\eta, \j) = p(\rho,\eta)+\frac{1}{2}\prn{\rho\alpha'(\rho)-\alpha(\rho)}\nrm{\j}^2, \quad  p(\rho,\eta) = \rho^2 \pd{\varepsilon}{\rho}.
\end{equation*} 
\subsection{Dissipation and asymptotic analysis}
It is possible to complement this system with \textit{dissipation}, at the level of the heat transfer. This can be done while still conserving the hyperbolic structure by using \textit{Rayleigh's dissipation} function \cite{goldstein2002classical}, so that the new system writes 
\begin{subequations}
	\begin{alignat}{2}
		& \frac{\partial \rho}{\partial t} &&+ \div{\rho \u } = 0, \\
		& \frac{\partial \rho\u}{\partial t}  &&+ \div{ \rho \u \otimes \u +\Pi  }= 0, \quad \Pi = P(\rho,\eta, \j) \ \mathbf{I}+\alpha(\rho)\ \j \otimes \j,  \\
		& \frac{\partial \j}{\partial t} &&+ \grad{\j\cdot\u + \theta(\rho,\eta)}  + \prn{ \frac{\partial \j}{\partial \x} - \prn{\frac{\partial \j}{\partial \x}}^T }\u = -\frac{\partial \calR}{\partial \j}, \label{dissipative_hyp_heat_eq_j} \\
		&\frac{\partial \rho\eta}{\partial t} &&+ \div{\rho \eta\u+ \alpha(\rho) \j  } = \frac{\alpha (\rho)}{\theta(\rho,\eta)}\,\frac{\partial \calR}{\partial\j}\cdot  \j. 
	\end{alignat}
	\label{eq:dissipative_hyp_heat_eq}
\end{subequations}
Here $\calR$ is the \textit{Rayleigh dissipation} function and which we take in the simplest form as
\begin{equation}
	\calR=\frac{1}{2\tau}\|\j\|^2, \qquad \pd{\calR}{\j} = \frac{1}{\tau}\j,
	\label{eq:Ray}
\end{equation}
where $\tau$ is a characteristic relaxation time. This structure is naturally compatible with the total energy conservation and straightforward computations, whose details are given in Appendix \ref{app:energy_conservation}, show that 
multiplying each of the equations in system \eqref{eq:dissipative_hyp_heat_eq} with the corresponding conjugate variable and summing up yields the energy conservation law
\begin{equation*}
	\pdt{E} + \div{E\u +\Pi \u+\mathbf{q}} = 0, \quad \mathbf{q} = \alpha(\rho)\, \theta(\rho,\eta)\,  \j,
\end{equation*}
where 
\begin{equation*}
	E=\frac{1}{2}\rho \nrm{\u}^2   +\rho\varepsilon(\rho,\eta)+\frac{1}{2}\alpha(\rho)\nrm{\j}^2.
\end{equation*}
The term $\mathbf{q}$ in the energy flux corresponds to interstitial working \cite{dunn1986thermomechanics}, and can be interpreted here as the heat flux density vector. In fact, the relaxation time $\tau$ can be chosen such that the energy equation is asymptotically consistent with Fourier's law of heat conduction. Assuming that $\tau$ is sufficiently small, we can expand all the variables in power series up to first order in $\tau$  as follows 
\begin{equation*}
	\rho = \rho_0 + O(\tau), \quad \u = \u_0 + O(\tau), \quad \eta = \eta_0 + O(\tau), \quad \j = \j_0 + \tau \j_1 + O(\tau^2) ,
\end{equation*}   
and then insert these expansions in equation \eqref{dissipative_hyp_heat_eq_j}. Multiplying the latter by $\tau$ allows us to write
\begin{equation*}
	\tau\prn{\pdt{\j_0}+\frac{\partial\j_0}{\partial \x} \u_0+ \prn{\frac{\partial \u_0}{\partial \x}}^T\j_0 + \nabla \theta(\rho_0,\eta_0)}   = -( \j_0 + \tau \j_1) + O(\tau^2).
\end{equation*}
Now, by equating in succession the coefficients of this expansion to zero, starting from the leading order, one obtains 
\begin{equation*}
	\j_0=0, \quad \j_1 = -\nabla \theta(\rho_0,\eta_0) \quad \Rightarrow \quad \j = -\tau\ \nabla \theta(\rho_0,\eta_0) + O(\tau^2).
\end{equation*} 
Under these considerations, the heat flux vector expresses as
\begin{equation*}
	\mathbf{q} = -\tau \alpha(\rho_0)\, \theta(\rho_0,\eta_0)\,  \nabla\theta(\rho_0,\eta_0).
\end{equation*}
Thus, we identify the latter with the one that appears in Fourier's law of heat conduction and which is given by $\mathbf{q}_F = -K \nabla \theta$ where $K$ is the thermal conductivity. As a result, one obtains an expression of the relaxation time
\begin{equation}
	\tau = \frac{K}{\alpha(\rho_0) \ \theta(\rho_0,\eta_0)},
	\label{eq:tau}
\end{equation}
for which our model is asymptotically consistent with Fourier's law of heat conduction. Let us remark that if we take $\tau$ as a constant, the constraint $\nabla\wedge \j=0$ still holds even in the presence of the right hand-side, provided the initial condition for $\j$ is also curl-free. Indeed applying the curl operator to equation \eqref{dissipative_hyp_heat_eq_j} in the case where $\calR$ is given by \eqref{eq:Ray} with $\tau=constant$ we get 
\begin{equation*}
	\pd{}{t}\prn{\nabla \wedge \j} + \nabla \wedge \prn{\prn{\nabla \wedge\j}\wedge \u} = -\frac{1}{\tau} \nabla \wedge \j.
\end{equation*}
If $\nabla \wedge \j\equiv0$ initially, it will remain as such for all $t>0$. In the one dimensional case, since the curl-free constraint is trivially satisfied, $\tau$ can also be a function of dependent variables.  

\subsection{Total energy convexity}
System of equations \eqref{eq:dissipative_hyp_heat_eq} is not a set of conservation laws since the evolution equation for $\j$ \eqref{dissipative_hyp_heat_eq_j} is not in conservative form. Therefore, proving hyperbolicity is not straightforward since we cannot use the Godunov-Lax-Friedrichs approach based on the existence of an additional scalar conservation law for a conservative system  \cite{godunov1961interesting,friedrichs1971systems}. It appears that the convexity of the energy as a function of $\prn{\rho\u,\j,\rho \eta,\rho}$ is nevertheless a sufficient condition of hyperbolicity (for proof, see Appendix \ref{app:sym}).  
Therefore, let us check the positive definiteness of the total energy hessian. The convexity of the total energy per unit volume is related to the convexity of the specific energy via the following theorem allowing to fairly simplify the algebra  \cite{godunov2003elements,wagner2009symmetric}.
\begin{thm} 
	\label{thm}
	Let $f(q_1,\cdots,q_n)$ be a convex function on a subset $\Omega\subset\mathbb{R_+}\times\mathbb{R}^{n-1}$. Then
	\begin{equation*}
		h\prn{1/q_1,q_2/q_1,\cdots,q_n/q_1} = 1/q_1\ f(q_1,q_2,\cdots,q_n)
	\end{equation*}
	is also convex.
\end{thm}
In our case, the specific energy considered as a function of $\prn{\u,v\j,\eta,v}$ is given by 
\begin{equation*}
	e(\u,v\j,\eta,v) = \frac{1}{2} \nrm{\u}^2  + \frac{\alpha(1/v)}{2v}\nrm{v \j}^2 +\varepsilon(v,\eta).
\end{equation*}
Here, by abuse of notation, we denoted $\varepsilon(1/v,\eta)$ by $\varepsilon(v,\eta)$. Furthermore, it suffices to take $\u=(u,0,0)$ and $\j=(j,0,0)$ as the considered energy only depends on the norms of these vectors. Under these considerations, the Hessian of $e$ in the basis $\prn{u,v j,\eta,v}$ is given by 
\begin{equation*}
	\everymath{\displaystyle}
	\nabla^2 e = \left(
	\begin{array}{cccc}
		1 & 0 & 0 & 0 \\
		0 & \varepsilon_{\eta\eta} & 0 & \varepsilon_{v\eta} \\
		0 & 0 & \alpha (1/v)/v & -j \left(\alpha '\left(1/v \right)+v  \alpha \left(1/v\right)\right)/v ^2 \\
		0 & \varepsilon_{v\eta} & -j \left(\alpha '\left(1/v \right)+v  \alpha \left(1/v\right)\right)/v ^2 & \varepsilon_{vv}+     \frac{j^2}{2v^3}   \prn{ \alpha ''\left(1/v\right)+4 v\alpha '\left(1/v\right)+ 2v^2 \alpha \left(1/v \right)} \\
	\end{array}
	\right).
\end{equation*}
In order to evaluate the positive definiteness of this matrix, we use here Sylvester's criterion, \textit{i.e.}, the determinant of the principal minors must be positive which yields the following inequalities
\begin{equation*}
	\varepsilon_{\eta\eta}>0, \qquad \alpha(1/v)/v>0, \qquad
	\prn{ \alpha (1/v) \alpha ''(1/v) - 2 \alpha'(1/v)^2 }    \frac{q^2}{2v^3}\varepsilon_{\eta\eta} + \alpha(1/v) \mathrm{det}\left(\nabla^2_{(v,\eta)}\varepsilon(v,\eta)\right)>0.
\end{equation*}
Taking into account that by definition $v>0$ and $\alpha(1/v)>0$ and assuming that $\varepsilon(v,\eta)$ is a convex function in $(v,\eta)$, it is therefore sufficient to have 
\begin{equation*}
	\alpha (1/v) \alpha ''(1/v) - 2 \alpha'(1/v)^2 \geq0 \quad \text{for }  v>0,
\end{equation*} 
which can be cast in an equivalent simpler form in terms of $\rho$
\begin{equation*}
	\pd{^2}{\rho^2}\prn{\frac{1}{\alpha(\rho)}} \leq 0, \quad \text{for } \rho > 0.
\end{equation*}
This simply means, that the energy $e$ is convex if $1/\alpha(\rho)$ is a concave function of $\rho$. 
One simple example verifying this criteria is $\alpha(\rho) = \varkappa^2/\rho$ where $\varkappa$ is a constant and which we will we consider in later parts of the paper.
\begin{rmrk*}
	As a consequence of theorem \ref{thm},  if $1/\alpha(\rho)$ is a concave function of $\rho$, then the energy 
	\begin{equation*}
		E(\mathbf{Q})= \frac{1}{2\rho} \nrm{\mathbf{m}}^2   +\rho\varepsilon(\rho,s/\rho)+\frac{1}{2}\alpha(\rho)\nrm{\j}^2
	\end{equation*}
	is also a convex function of $\mathbf{Q}^T = \prn{\rho,\mathbf{m}^T,s,\j^T}$ where $\mathbf{m}=\rho \u$ and $s=\rho \eta$.
\end{rmrk*}

\subsection{Friedrichs symmetric form}
Symmetric $t-$hyperbolic systems in the sense of Friedrichs  \cite{friedrichs1954symmetric} is a very important class of partial differential equations. Not only it allows to establish the well-posedness of the Cauchy problem \cite{kato1975cauchy}, it also simplifies linear wave analysis in terms of symmetric matrices. Let us introduce the conjugate variables to $\mathbf{Q}$ defined by 
\begin{equation*}
	\mathbf{U}^T = \pd{E}{\mathbf{Q}} = \prn{E_\rho,E_{\mathbf{m}},E_{s},E_\j} =  \prn{\mu, \u^T,\theta, \alpha\,\j^T}.
\end{equation*}
where $\mu=E_\rho$ is the generalized chemical potential.
Let $E^\star$ be the Legendre transform of the energy $E$
\begin{equation*}
	E^\star(\mathbf{U})= \mathbf{U}^T\mathbf{Q}-E(\mathbf{Q}), \qquad \mathbf{Q}^T = \pd{E^\star}{\mathbf{U}}. 
\end{equation*}
In the variables $\mathbf{U}$, the system is symmetric hyperbolic in the sense of Friedrichs if $E(\mathbf{Q})$ is a convex function of $\mathbf{Q}$. Under these notations, system of equations \eqref{eq:dissipative_hyp_heat_eq} can be cast into the symmetric form (see Appendix \ref{app:sym}) 
\begin{equation*}
	\tilde{\mathbf{A}}\pdt{\mathbf{U}} + \sum_i \tilde{\mathbf{B}}^i \pd{\mathbf{U}}{x_i} = \mathbf{S}, 
\end{equation*}
where $\tilde{\mathbf{A}}$ is the Hessian matrix of $E^\star(\mathbf{U})$, $\tilde{\mathbf{B}}^i$ are symmetric matrices related to the fluxes and non-conservative terms in the governing equations and $\mathbf{S}$ is the vector of source terms. 

\section{Hyperbolic model analysis}
\label{sec:Hyperbolic}
\subsection{Eigenstructure study} 
System of equations \eqref{eq:dissipative_hyp_heat_eq} is provably invariant by rotations of the $\mathcal{SO}(3)$ group. This allows us to reduce the study of its hyperbolicity to the one-dimensional case, \textit{i.e.}, all the components of all the vector variables are kept in the three dimensions of space but the variables are functions of only $(x,t)$ instead of $(x,y,z,t)$. It suffices to study the homogeneous system, written here in quasilinear form
\begin{equation*}
	\pdt{\mathbf{V}} + \mathbf{A}(\mathbf{V}) \pd{\mathbf{V}}{\x} = 0,
\end{equation*}
where $\mathbf{V}$ is the vector of primitive variables and $\mathbf{A}$ is the quasilinear matrix, both given by
\begin{equation*}
	\mathbf{V} = \prn{\begin{array}{c}
			\rho\\
			u_1\\
			u_2\\
			u_3\\
			\eta\\
			j_1\\
			j_2\\
			j_3
	\end{array}}, \quad
	\mathbf{A} = \left(
	\begin{array}{cccccccc}
		u_1 & \rho  & 0 & 0 & 0 & 0 & 0 & 0 \\
		\frac{1}{\rho}\pd{p}{\rho}+\varkappa^2\frac{\prn{j_2^2+j_3^2}}{\rho^3} & u_1 & 0 & 0 & \frac{1}{\rho}\pd{p}{\eta} & 0 & -2 \varkappa^2\frac{  j_2}{\rho ^2} & - 2 \varkappa^2 \frac{ j_3}{\rho ^2} \\
		-\varkappa^2  \frac{j_1 j_2}{\rho ^3} & 0 & u_1 & 0 & 0 & \varkappa^2 \frac{  j_2}{\rho ^2} & \varkappa^2 \frac{  j_1}{\rho ^2} & 0 \\
		-\varkappa^2 \frac{  j_1 j_3}{\rho ^3} & 0 & 0 & u_1 & 0 & \varkappa^2 \frac{  j_3}{\rho ^2} & 0 & \varkappa^2 \frac{  j_1}{\rho ^2} \\
		-\varkappa^2 \frac{  j_1}{\rho ^3} & 0 & 0 & 0 & u_1 & \varkappa^2 \frac{1 }{\rho ^2} & 0 & 0 \\
		\pd{\theta}{\rho} & j_1 & j_2 & j_3 & \pd{\theta}{\eta} & u_1 & 0 & 0 \\
		0 & 0 & 0 & 0 & 0 & 0 & u_1 & 0 \\
		0 & 0 & 0 & 0 & 0 & 0 & 0 & u_1 \\
	\end{array}
	\right).
\end{equation*}
Straightforward linear algebra shows that $\mathbf{A}$ admits 8 eigenvalues whose expressions are given by  
\begin{align*}
	\everymath{\displaystyle}
	\begin{cases}
		\chi_1 = u_1 - \sqrt{Z_1+Z_2},\\
		\chi_2 = u_1 - \sqrt{Z_1-Z_2},\\
		\chi_{3-6} = u_1,\\
		\chi _7 = u_1 + \sqrt{Z_1-Z_2},\\
		\chi _8 = u_1 + \sqrt{Z_1+Z_2}
	\end{cases}
	\quad \text{where} \quad 
	\begin{cases}
		Z_1 = \frac{1}{2}\prn{a_p^2 + a_T^2+a_q^2},  \\[2.5mm]
		Z_3 = \frac{1}{2}\prn{a_p^2 - a_T^2}, \quad Z_2 = \sqrt{ a_{pT}^4+Z_3^2}, \\[2.5mm]
		a_p^2 = \pd{p}{\rho}=v^2\varepsilon_{vv}, \quad a_T^2 = \frac{\varkappa^2}{\rho^2}  \pd{\theta}{\eta} = \varkappa^2v^2\varepsilon_{\eta\eta},\\[2.5mm]
		a_{pT}^4= \frac{\varkappa^2}{\rho^2} \pd{p}{\eta}\pd{\theta}{\rho} = \varkappa^2v^4\varepsilon_{v\eta}^2  , \quad 	a_q^2 =  \frac{2\varkappa^2}{\rho^2}\prn{j_2^2+j_3^2},
	\end{cases}
\end{align*}
where $a_p$, $a_T$, $a_{pT}$ and $a_q$ are auxiliary quantities homogeneous to velocities.
In order to prove hyperbolicity for this system, we need first to prove that all the eigenvalues $\chi_{1-8}$ are real-valued. Given that both $Z_1$ and $Z_2$ are trivially positive, it suffices to prove in this case that $Z_1-Z_2\geq0$. This follows immediately by rewriting 
\begin{equation*}
	Z_2 = \sqrt{Z_1^2-\prn{a_p^2a_T^2-a_{pT}^4}},
\end{equation*}
and by identifying that $a_p^2a_T^2-a_{pT}^4 = \varkappa^2v^4\mathrm{det}\prn{\nabla^2_{(v,\eta)}\varepsilon(v,\eta)}\geq0$, thus leading to $Z_1\geq Z_2$. Therefore, the eigenvalues $\chi_{1-8}$ are real. It remains to check, in the general case whether $\mathbf{A}$ admits a full corresponding basis of right eigenvectors. This is not the case, as straightforward computations show that $\mathbf{A}$ is missing two right eigenvectors for the multiple eigenvalue $u_1$, meaning that system \eqref{eq:dissipative_hyp_heat_eq} is only weakly hyperbolic. This does not come out as a surprise since this shortcoming appears in several models where an independent field bound by a curl-free constraint is evolved  \cite{SHTCSurfaceTension,Schmidmayer2017,dhaouadi2022NSK,dhaouadiThinFilms,romenski2010conservative,chiocchetti2023exactly,riomartin2023highorder}. Nonetheless, there are procedures allowing to recover strong hyperbolicity in this case. One can for instance apply a curl-cleaning procedure, by introducing an artificial vector field $\ppsi$ to the equations allowing to propagate discrete curl-errors to domain boundaries \cite{MunzCleaning}. In this case, the full system of equations accounting for the curl-cleaning is given below   
\begin{subequations}
	\begin{alignat}{2}
		& \frac{\partial \rho}{\partial t} &&+ \div{\rho \u } = 0, \\
		& \frac{\partial \rho\u}{\partial t}  &&+ \div{ \rho \u \otimes \u + P(\rho,\j) \ \mathbf{I}+\frac{\varkappa^2}{\rho}\ \j \otimes \j }= 0, \\
		&\frac{\partial \rho\eta}{\partial t} &&+ \div{\rho \eta\u+ \frac{\varkappa^2}{\rho} \j  } = 0,\\
		& \frac{\partial \j}{\partial t} &&+ \grad{\j\cdot\u + \theta(\rho,\eta)}  + \prn{ \frac{\partial \j}{\partial \x} - \prn{\frac{\partial \j}{\partial \x}}^T }\u + \frac{\rho }{\varkappa}a_c \nabla\wedge\ppsi= 0, \label{eq:p_clean} \\
		& \frac{\partial \ppsi}{\partial t} && + \pd{\ppsi}{\x}\u  - \frac{\varkappa}{\rho}a_c \nabla \wedge \j= 0. \label{eq:psi_clean}
	\end{alignat}
	\label{eq:heat_eq_+clean}
\end{subequations}
where $a_c$ is the \textit{cleaning velocity}. The factors multiplying the curl terms in equations \eqref{eq:p_clean} and \eqref{eq:psi_clean} are not arbitrary. They are chosen as such in order to ensure thermodynamic compatibility of this modified system of equations with total energy conservation, which in this case becomes (see Appendix \ref{app:energy_conservation})
\begin{equation*}
	\pdt{\calE} + \div{\calE\u +\Pi \u+\frac{\varkappa^2}{\rho} \theta(\rho,\eta)\,  \j + \varkappa a_c \j \wedge \ppsi} = 0,
\end{equation*}
where the modified energy writes as 
\begin{equation}
	\calE = \frac{1}{2}\rho\nrm{\mathbf{u}}^2 + \frac{\varkappa^2}{2\rho}\nrm{\j}^2 + \rho\varepsilon(\rho,\eta) + \frac{1}{2}\rho \nrm{\ppsi}^2,
	\label{eq:energy_cleaning}
\end{equation}
In three space dimensions, this system admits the eigenvalues $\chi^c_i$ which are given by 
\begin{gather*}
	\chi^c_{1,2} = u_1 - a_c,\qquad \chi^c_3 = u_1 - \sqrt{Z_1+Z_2}, \qquad
	\chi^c_4 = u_1 - \sqrt{Z_1-Z_2},\qquad
	\chi^c_{5-7} = u_1, \\
	\chi^c_8 = u_1 + \sqrt{Z_1-Z_2}, \qquad
	\chi^c_9 = u_1 + \sqrt{Z_1+Z_2}, \qquad
	\chi^c_{10,11}= u_1 +a_c. 
\end{gather*} 
To these eigenvalues corresponds a full basis of right eigenvectors and which are collected below in the corresponding order of the eigenvalues
\begin{equation*}
	\small
	\mathbf{R}^c=\left(
	\begin{array}{ccccccccccc}
		j_3 r_a^+ & -j_2r_a^+ & 1 & 1 & 0 & \frac{\rho }{j_1} & 0 &1 & 1 & j_2 r_a^- & -j_3r_a^- \\
		-\frac{j_3 a_c r_a^+}{\rho } & \frac{j_2 a_c r_a^+}{\rho } & -\frac{\varphi _{12}^+ }{\rho } & -\frac{\varphi _{12}^- }{\rho } & 0 & 0 & 0 & \frac{\varphi _{12}^- }{\rho } & \frac{\varphi _{12}^+ }{\rho } & \frac{j_2 a_c r_a^-}{\rho } & -\frac{j_3 a_c r_a^-}{\rho } \\
		\frac{\varkappa  j_2 j_3 r_b^+}{\rho ^2} & -\frac{\varkappa  j_1}{\rho  a_c} - \frac{\varkappa  j_2^2 r_b^+}{\rho ^2} & \frac{\varkappa  j_2 }{\rho ^2}\frac{\varphi _{23}^-}{\varphi _{23}^+} & -\frac{ \varkappa  j_2 }{\rho ^2} \frac{\varphi _{23}^+}{\varphi _{23}^-} & 0 & r_{36} & -\frac{j_3}{j_2} & -\frac{\varkappa  j_2 }{\rho ^2} \frac{\varphi _{23}^+}{\varphi _{23}^-}&  \frac{ \varkappa  j_2 }{\rho ^2}\frac{\varphi _{23}^-}{\varphi _{23}^+} & -\frac{\varkappa  j_1}{\rho  a_c}-\frac{\varkappa  j_2^2 r_b^-}{\rho ^2} & \frac{\varkappa  j_2 j_3 r_b^-}{\rho ^2} \\
		\frac{\varkappa  j_1}{\rho  a_c}+\frac{\varkappa  j_3^2 r_b^+}{\rho ^2} & -\frac{\varkappa  j_2 j_3 r_b^+}{\rho ^2} & \frac{ \varkappa  j_3 }{\rho ^2}\frac{\varphi _{23}^-}{\varphi _{23}^+} & -\frac{ \varkappa  j_3 }{\rho ^2}\frac{\varphi _{23}^+}{\varphi _{23}^-} & 0 & 0 & 1 & -\frac{ \varkappa  j_3 }{\rho ^2}\frac{\varphi _{23}^+}{\varphi _{23}^-} & \frac{ \varkappa  j_3 }{\rho ^2} \frac{\varphi _{23}^-}{\varphi _{23}^+}& -\frac{\varkappa  j_2 j_3 r_b^-}{\rho ^2} & \frac{\varkappa  j_1}{\rho  a_c}+\frac{\varkappa  j_3^2 r_b^-}{\rho ^2} \\
		\frac{\varkappa  j_3 r_b^+}{\rho ^2} & -\frac{\varkappa  j_2 r_b^+}{\rho ^2} & \frac{ \varkappa  }{\rho ^2} \frac{\varphi _{23}^-}{\varphi _{23}^+}& -\frac{\varkappa  }{\rho ^2}\frac{\varphi _{23}^+}{\varphi _{23}^-}  & 0 & r_{56} & 0 & -\frac{ \varkappa  }{\rho ^2}\frac{\varphi _{23}^+}{\varphi _{23}^-} & \frac{\varkappa  }{\rho ^2}\frac{\varphi _{23}^-}{\varphi _{23}^+}  & -\frac{\varkappa  j_2 r_b^-}{\rho ^2} & \frac{\varkappa  j_3 r_b^-}{\rho ^2} \\
		\frac{j_1 j_3 r_a^+}{\rho }-\frac{j_3 a_c r_b^+}{\varkappa } & \frac{j_2 a_c r_b^+}{\varkappa }-\frac{j_1 j_2 r_a^+}{\rho } & r^{-}_+ & r^{+}_- & 0 & 1 & 0 & r^{-}_- & r^{+}_+ & \frac{j_1 j_2 r_a^-}{\rho }-\frac{j_2 a_c r_b^-}{\varkappa } & \frac{j_3 a_c r_b^-}{\varkappa }-\frac{j_1 j_3 r_a^-}{\rho } \\
		0 & \frac{\rho }{\varkappa } & 0 & 0 & 0 & 0 & 0 & 0 & 0 & -\frac{\rho }{\varkappa } & 0 \\
		-\frac{\rho }{\varkappa } & 0 & 0 & 0 & 0 & 0 & 0 & 0 & 0 & 0 & \frac{\rho }{\varkappa } \\
		0 & 0 & 0 & 0 & 1 & 0 & 0 & 0 & 0 & 0 & 0 \\
		1 & 0 & 0 & 0 & 0 & 0 & 0 & 0 & 0 & 0 & 1 \\
		0 & 1 & 0 & 0 & 0 & 0 & 0 & 0 & 0 & 1 & 0 \\
	\end{array}
	\right).
\end{equation*}

where auxiliary quantities where introduced to lighten the expressions and whose definitions are as follows
\begin{gather*}
	\varphi_{ij}^{\pm}=\sqrt{\phi_i\pm \phi_j},\quad r_{+}^{\pm} =  \frac{j_1}{\rho}\pm  \frac{\varphi^-_{12} \varphi^+_{23}}{\varkappa \ \varphi^-_{23}}, \quad 	r_{-}^{\pm} = \frac{j_1}{\rho}\pm\frac{\varphi^+_{12} \varphi^-_{23}}{\varkappa \ \varphi^+_{23}}, 	\quad r_{36} = \rho\frac{  \prn{\varphi^+_{12}\varphi^-_{12}}^2-\frac{1}{2}a_q^2 \prn{\varphi^+_{13}}^2}{\varkappa    j_1 j_2 a_{pT}}, \quad r_{56} = -\frac{\varkappa  \prn{\varphi^+_{13}}^2}{\rho  j_1 \varphi^-_{23}\varphi^+_{23}}, \\
	r_a^\pm = \frac{\varkappa}{\rho a_c}\frac{  2 a_c \rho  \left(a_c^2-\prn{\varphi^-_{13}}^2\right)\pm\varkappa  j_1 \varphi^-_{23}\varphi^+_{23}}{\prn{a_c^2-\prn{\varphi_{12}^+}^2}\prn{a_c^2-\prn{\varphi_{12}^-}^2}}, \quad r_b^\pm = \frac{\varkappa}{\rho a_c}\frac{\varkappa  j_1 \left(a_c^2-\prn{\varphi^+_{13}}^2\right)\pm2 a_c \rho  \varphi^-_{23}\varphi^+_{23}}{\prn{a_c^2-\prn{\varphi_{12}^+}^2}\prn{a_c^2-\prn{\varphi_{12}^-}^2}}.
\end{gather*}

\subsection{Structure of the eigenfields}
\label{subs:nature}
We consider here the homogeneous system \eqref{eq:hyp_heat_eq} in one dimension of space and which writes
	\begin{alignat*}{2}
		& \frac{\partial \rho}{\partial t} &&+ u\pd{\rho }{x} + \rho\pd{u}{x} = 0, \\
		& \frac{\partial u}{\partial t}  &&+ u\pd{u }{x} + \frac{1}{\rho} \pd{p}{\rho}\pd{\rho}{x} + \frac{1}{\rho} \pd{p}{\eta}\pd{\eta}{x}= 0 , \\
		&\frac{\partial \eta}{\partial t} &&+ u\pd{\eta}{x}+\frac{\varkappa^2}{\rho^2}\pd{j}{x} -\frac{\varkappa ^2 }{\rho ^3}j \pd{\rho}{x}  = 0, \\
		& \frac{\partial j}{\partial t} &&+ j\pd{u}{x} + u\pd{j}{x} +  \pd{\theta}{\rho}\pd{\rho}{x} + \pd{\theta}{\eta}\pd{\eta}{x}= 0. 
	\end{alignat*}
The associated quasilinear matrix is given by 
\begin{align*}
	\everymath{\displaystyle}
	\mathbf{A}_{1d}=\left(
	\begin{array}{cccc}
		u & \rho  & 0 & 0 \\[2mm]
		\frac{1}{\rho}\pd{p}{\rho} & u & \frac{1}{\rho}\pd{p}{\eta} & 0 \\[2mm]
		-\frac{\varkappa ^2 }{\rho ^3}j & 0 & u & \frac{\varkappa ^2}{\rho ^2} \\[2mm]
		\pd{\theta}{\rho} & j & \pd{\theta}{\eta} & u
	\end{array}
	\right),
\end{align*}
and the corresponding eigenvalues, presented in the same fashion as previously, are 
\begin{equation*} 
	\begin{cases}
		\lambda_1 = u - \sqrt{Y_1+Y_2},\\
		\lambda_2 = u - \sqrt{Y_1-Y_2},\\
		\lambda _3 = u + \sqrt{Y_1-Y_2},\\
		\lambda _4 = u + \sqrt{Y_1+Y_2},
	\end{cases}
	\quad \text{where} \quad 
	\begin{cases}
		Y_1 = \frac{1}{2}\prn{a_p^2 + a_T^2}, \\
		Y_2 = \sqrt{ a_{pT}^4+Y_3^2}, \\
		Y_3 = \frac{1}{2}\prn{a_p^2 - a_T^2}.
	\end{cases}
\end{equation*}
The eigenvalues $\lambda_{1-4}$ are distinct for $\varkappa\neq0$. They can be categorized into faster characteristics $\lambda_{1,4}$ and slower characteristics $\lambda_{2,3}$, corresponding respectively to the propagation of \textit{acoustic} and \textit{thermal} waves. In the same context, we will refer to the quantities $\sqrt{Y_1+Y_2}$ and $\sqrt{Y_1-Y_2}$ as the \textit{acoustic} and \textit{thermal} sound speeds, respectively.
One can compute, the associated right eigenvectors which are collected hereafter as the columns of the following matrix, from the left to the right in the same order as their corresponding eigenvalues
\begin{equation*}
	\everymath{\displaystyle}
	\mathbf{r} = \left(
	\arraycolsep=4pt
	\begin{array}{cccc}
		\rho  & \rho  & \rho  & \rho  \\[1.5mm]
		-\upsilon_{12}^+ & -\upsilon_{12}^- & \upsilon_{12}^- & \upsilon_{12}^+ \\[2mm]
		\frac{\varkappa}{\rho}\frac{  \upsilon_{23}^-}{  \upsilon_{23}^+} & -\frac{\varkappa}{\rho}\frac{ \upsilon_{23}^+}{  \upsilon_{23}^-} & -\frac{\varkappa}{\rho}\frac{ \upsilon_{23}^+}{ \upsilon_{23}^-} & \frac{\varkappa}{\rho}\frac{  \upsilon_{23}^-}{ \upsilon_{23}^+} \\[4mm]
		j-\frac{\rho}{\varkappa}\frac{\upsilon_{12}^+ \upsilon_{23}^-}{ \upsilon_{23}^+} & j+\frac{\rho}{\varkappa}\frac{\upsilon_{12}^- \upsilon_{23}^+}{  \upsilon_{23}^-} & j-\frac{\rho}{\varkappa}\frac{  \upsilon_{12}^- \upsilon_{23}^+}{ \upsilon_{23}^-} & j+\frac{\rho}{\varkappa}\frac{ \upsilon_{12}^+ \upsilon_{23}^-}{  \upsilon_{23}^+} \\
	\end{array}
	\right),
\end{equation*}
where $\upsilon_{mn}^{\pm} = \sqrt{Y_m\pm Y_n}$. As an example, let us consider the case of a polytropic gas, for which the internal energy is 
\begin{equation*}
	\varepsilon(\rho,\eta) = \frac{\rho^{\gamma-1}}{\gamma-1}e^{\eta/c_V},
\end{equation*}
where $c_V$ is the specific heat capacity at constant volume and $\gamma$ is the heat capacity ratio. One can check in this case that the first and fourth eigenfields are genuinely nonlinear in the sense Lax  \cite{lax1973hyperbolic} while the second and the third eigenfields are neither genuinely nonlinear not linearly degenerate. Let us illustrate the latter property on the particular case $\gamma=2$ and $c_V=1$. Direct computations show that while $\nabla_\mathbf{v}\lambda_4\cdot \mathbf{r}_4 \neq0$ for all admissible states, $\nabla_\mathbf{v}\lambda_3\cdot \mathbf{r}_3$ is equal to zero on the hyperplane of exact equation
\begin{equation}
	\varkappa\;v = \frac{1}{\sqrt{6}}\sqrt{2+\sqrt[3]{17-12 \sqrt{2}}+\sqrt[3]{17+12 \sqrt{2}}}= 1.0399...,
	\label{eq:hyperplane}
\end{equation}
which separates two regions of different signs of $\nabla_\mathbf{v}\lambda_3\cdot \mathbf{r}_3$. This yields the preliminary conclusion that expansion shocks and compression fans corresponding to the thermal characteristic fields can exist as solutions to the Riemann problem for our system  \cite{zel1946,oleinik1959uniqueness,liu1981admissible,borisov1983rarefaction,menikoff1989riemann,lefloch2002hyperbolic,kluwick2018shock}. This will be further developed in section \ref{sec:shocks}.

\subsection{Dispersion relation}
\begin{figure}[H]
	\includegraphics[width=0.5\textwidth]{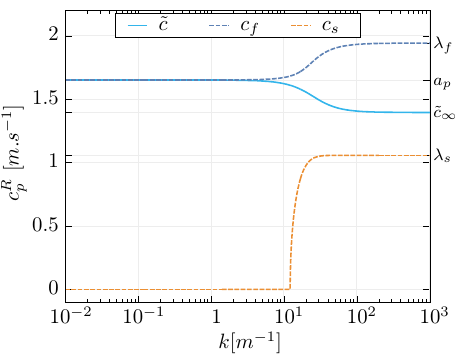}
	\includegraphics[width=0.5\textwidth]{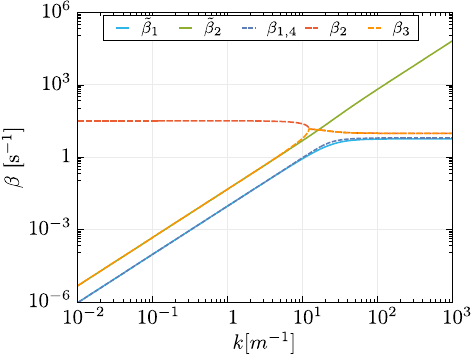}
	\caption{log-linear plot of the real part of the phase velocities (Left) and log-log plot of the attenuation factors (Right) for both system of equations \eqref{eq:dissipative_hyp_heat_eq} (Dashed lines) and system \eqref{eq:Euler-Fourier} (Solid lines) as a function of the wavenumber $k$. Polytropic gas equation of state is used here with $\gamma=1.4$, $c_V=3/2$, $\rho_0=1$, $\eta_0=1$, $\varkappa=1$ and $K=0.1$.}
	\label{fig:cp}
\end{figure}

In order to further investigate the properties of the proposed model and also assess its stability, we study hereafter its dispersion relation. 
System of equations \eqref{eq:dissipative_hyp_heat_eq} is provably invariant by rotations of the $\mathcal{SO}(3)$ group. This allows us to reduce the analysis to the one-dimensional case.
Let us consider a reference equilibrium state  defined by $\mathbf{V}_0=(\rho_0,u_0,j_0,\eta_0)^T$, with $u_0=0, j_0=0$. For convenience, we shift here to the corresponding set of conjugate variables $\mathbf{U}_0$, as it allows for a more straightforward verification of stability. Let $\mathbf{U}'$ be a small perturbation around this equilibrium. We are looking for the solution  $\mathbf{U}$ in the form :
\begin{equation*}
	\mathbf{U} = \mathbf{U}_0 + \mathbf{U}'.
\end{equation*}
Then, we linearize the one-dimensional system around $\mathbf{U}_0$ to obtain 
\begin{equation*}
	\tilde{\mathbf{A}}_{0}	\pdt{\mathbf{U}'} + \tilde{\mathbf{B}}_{0}\pd{\mathbf{U}'}{x} = \tilde{\mathbf{C}}_0 \mathbf{U}', \quad \tilde{\mathbf{A}}_{0}^T=\tilde{\mathbf{A}}_{0}>0,\quad \tilde{\mathbf{B}}_{0}^T=\tilde{\mathbf{B}}_{0}, \quad \tilde{\mathbf{C}}_0^T=\tilde{\mathbf{C}}_0\le 0, 
\end{equation*}
where the subscript $0$ means that the matrices are evaluated at the equilibrium state $\mathbf{U}_0$. Then we look for harmonic wave solutions of the form
\begin{equation*}
	\mathbf{U}' = \mathbf{U}_1 e^{i(kx-\omega t)}, \quad \text{with}\ k\in \mathbb{R}_+,\ \omega \in \mathbb{C},
\end{equation*}
which allows us to obtain the algebraic equation
\begin{equation}
	\prn{-\frac{\omega}{k} 	\tilde{\mathbf{A}}_{0} +\tilde{\mathbf{B}}_0+\frac{i}{k}\tilde{\mathbf{C}}_0}\mathbf{U}_1=0.
	\label{eq:disp_algeb}
\end{equation}
We define the phase velocity as the following complex quantity
\begin{equation*}
	c_p = \frac{\omega}{k}= c_p^R + i \ c_p^I, \quad \text{with} \ c_p^R,c_p^I \in \mathbb{R}.
\end{equation*}
In this setting, the stability of the solutions naturally requires $c_p^I\leq0$. This can be proven as follows, similarly to what was proposed in \cite{gavrilyuk2005acoustic}. We multiply equation \eqref{eq:disp_algeb} by $\overline{\mathbf{U}}^{T}$ to the left, where the bar denotes the complex conjugate. This allows to write  
\begin{equation*}
	\overline{\mathbf{U}}_1^{T} 	\prn{-c_p^R\tilde{\mathbf{A}}_{0} +\tilde{\mathbf{B}}_0 }\mathbf{U}_1 - i\prn{c_p^I\ \	\overline{\mathbf{U}}_1^{T} 	\tilde{\mathbf{A}}_{0}  \mathbf{U}_1 - \frac{1}{k}	\overline{\mathbf{U}}_1^{T}\tilde{\mathbf{C}}_0\mathbf{U}_1}=0.
\end{equation*}
Recall that in the conjugate variables $\mathbf{U}$, the matrices $\tilde{\mathbf{A}}_{0}$ and $\tilde{\mathbf{B}}_{0}$ are symmetric matrices meaning that the first term is real. This allows to obtain an expression of $c_p^I$ in the form 
\begin{equation*}
	c_p^I= \frac{1}{k}\frac{\mathcal{R}e\prn{\overline{\mathbf{U}}_1^{T}\tilde{\mathbf{C}}_0\mathbf{U}_1}}{\overline{\mathbf{U}}_1^{T} 	\tilde{\mathbf{A}}_{0}  \mathbf{U}_1} \le 0, \quad \text{since} \quad \tilde{\mathbf{C}}_0 \le 0.
\end{equation*}
The non-negative coefficient $\beta=-k c_p^I$ is also commonly referred to as the wave attenuation factor. Now we would like to evaluate both of $c_p^R$ and $c_p^I$ for our system and compare them with their equivalents for the Euler system with heat transfer. The phase velocities for our system, denoted by $c_{p_j}, j\in\{1..4\}$ can be obtained as the roots of the following polynomial 
\begin{equation*}
	F(c_p)=\left.\prn{c_p^4 +\frac{i }{k \tau }c_p^3 - \prn{\frac{\varkappa ^2}{ \rho^2  } \pd{\theta}{\eta} +  \pd{p}{\rho}}c_p^2  -\frac{i }{k \tau }  \pd{p}{\rho} c_p+\frac{\varkappa ^2}{\rho^2}  \prn{\pd{p}{\rho} \pd{\theta}{\eta}-\pd{p}{\eta} \pd{\theta}{\rho}}}\right\vert_{\mathbf{U}=\mathbf{U}_0},
\end{equation*}
while the phase velocities for the Euler system with heat transfer, linearized around $\tilde{\mathbf{U}}_0=\prn{{\rho_0,u_0,\eta_0}}$, denoted by $\hat{c}_{p_j}, j\in\{1..3\}$, are also obtainable from the following polynomial
\begin{equation*}
	\tilde{F}(\tilde{c}_p)=\left.\prn{\tilde{c}_p^3 + \frac{i K  k }{\rho  \theta} \pd{\theta}{\eta}\tilde{c}_p^2 -\pd{p}{\rho} \tilde{c}_p -\frac{i K k }{\rho  \theta} \prn{\pd{p}{\rho} \pd{\theta}{\eta}  -\pd{p}{\eta} \pd{\theta}{\rho}}}\right\vert_{\tilde{\mathbf{V}}=\tilde{\mathbf{V}}_0} .
\end{equation*}
It is worthy of note that for the value $\varkappa^2=K\rho_0\tau/\theta$ which gives asymptotic consistency of both models \eqref{eq:tau}, one has the algebraic relation
\begin{equation*}
	F(c_p)=c_p^4-\pd{p}{\rho}c_p^2 +\frac{i}{k \tau}\tilde{F}(c_p).
\end{equation*}
Recall that for a hyperbolic system with source terms, in the short wave limit $(k\rightarrow\infty)$, the real parts of phase velocities converge towards characteristic speeds. In our case, we have noted them $\lambda_{1-4}$ such that $\lambda_1<\lambda_2<0<\lambda_3<\lambda_4$. At the rest state we further have $\lambda_4=-\lambda_1=\lambda_f$ and $\lambda_3=-\lambda_2=\lambda_s$ where the subscripts stand for \textit{fast} and \textit{slow}, respectively. Therefore, it makes sense to also categorize the real parts of the phase velocities accordingly into fast and slow velocities and decompose them as follows
\begin{equation*}
	c_{p_1} = -c_f - i \beta_1/k, \quad c_{p_2} = -c_s - i \beta_2/k, \quad c_{p_3} = c_s - i \beta_3/k, \quad c_{p_4} = c_f - i \beta_4/k, 
\end{equation*}
where $c_s,\;c_f,\;\beta_{1-4}$ are non-negative quantities representing the slow velocities, fast velocities and attenuation factors respectively. In a similar way, we can also write the Euler system's eigenvalues as 
\begin{equation*}
	\tilde{c}_{p_{1}} = \tilde{c} -i\tilde\beta_1/k, \quad \tilde{c}_{p_{2}} = -\tilde{c} -i\tilde\beta_1 /k\quad \tilde{c}_{p_3} = -i\tilde\beta_2/k,
\end{equation*}
where $\tilde{c},\;\tilde{\beta}_{1-2}$ are non-negative.
Figure \ref{fig:cp} provides a graphical representation of the real parts $c_f,c_s$ and $\tilde{c}$ as well the attenuation factors $\beta_j$ and $ \tilde{\beta}_j$, as a function of the wavenumber $k$ for the particular case of a polytropic gas. 
As can be seen from Figure \ref{fig:cp}, the phase velocities of the proposed model are in agreement with the Euler system with heat conduction, both at the level of the wave speeds as well as the attenuation factors, at least in the long wave regime and up to a certain cut-off wavenumber. Note that the latter is an increasing function of $\varkappa$. After this threshold, the curves eventually separate since they have different limits when $k\rightarrow\infty$. Indeed, if we consider for instance the velocity $\tilde{c}$, to which corresponds $c_f$, the limit for $k\rightarrow0$ is the same for both and is the \textit{adiabatic} sound speed $a_p$
\begin{equation*}
	\lim\limits_{k\rightarrow 0}\tilde{c} = \lim\limits_{k\rightarrow 0}c_{f} = a_p.
\end{equation*}   
However, the limits for $k\rightarrow \infty$ are 
\begin{equation*}
	\lim\limits_{k\rightarrow \infty}\tilde{c} =\sqrt{\pd{p}{\rho}  -\pd{p}{\eta} \pd{\theta}{\rho}/\pd{\theta}{\eta}} = \sqrt{ \rho^{\gamma-1} e^{\eta/cv}} = \tilde{c}_\infty , \qquad  	\lim\limits_{k\rightarrow \infty}c_{f} = \lambda_f.
\end{equation*}   
%\begin{figure}[H]
%	\includegraphics[width=0.5\textwidth]{cp_kappa}
%	\includegraphics[width=0.5\textwidth]{beta_kappa}
%	\caption{On the left, log-log plot of $c_f$ (Solid), $\tilde{c}$(Dashed). On the right, log-log plot of $\beta_{1,3}$ (Solid), $\tilde{\beta}_{1,2}$ (Dashed). The plots are done for different values of $\varkappa$. Polytropic gas equation of state is used here with $\gamma=1.4$, $c_V=3/2$, $\rho_0=1$, $\eta_0=1$ and $K=0.1$.}
%	\label{fig:cp_kappa}
%\end{figure}

\section{Admissibility of shock waves and nonexistence of contact discontinuities}
\label{sec:shocks}
In this section we study the Rankine-Hugoniot conditions in an attempt to analyze the structure of shocks that are admissible to the system. In this context let us consider a discontinuity moving at constant speed $\mathcal{D}$, separating two different states. For physical relevancy, the conservation of mass, momentum, energy and $j$ must be satisfied. The entropy inequality will be considered as an admissibility criterion. We show below that in this case contact discontinuities do not exist for $\varkappa\neq0$. Note that since we are interested mainly in the jump relations, the relaxation sources will be omitted in the analysis. The jump of an arbitrary function $g$ will be denoted as $\RH{g} = g_R-g_L$ where $g_L$ and $g_R$ are the states on the left and on the right of the discontinuity, respectively. For the chosen set of conserved variables, the Rankine-Hugoniot conditions write 
\begin{subequations}
	\begin{align}
		\RH{\calM} &= 0, \label{eq:RH_rho} \\
		\RH{p+\frac{\calM^2}{\rho}} &= 0, \label{eq:RH_rhou}\\
		\RH{\calM\prn{\frac{\calM^2}{2\rho^2} + \varepsilon+\frac{p}{\rho} + \frac{1}{2}\frac{\varkappa^2}{\rho^2}j^2} + \frac{\varkappa^2}{\rho}\theta\; j }&= 0, \label{eq:RH_E}\\
		\RH{\calM \frac{j}{\rho}+\theta} &=0,	\label{eq:RH_j}
	\end{align}
	\label{eq:RH}
\end{subequations}
where we have defined the mass flux across the discontinuity front by $\calM = \rho (u-\calD)$. 
\subsection{Non-existence of contact discontinuities}
On contact discontinuities, that is for $\calM=0$, one obtains by direct substitution
\begin{equation*}
	\RH{p}=0, \quad \RH{\varkappa^2\frac{j}{\rho}\theta}=0, \quad \RH{\theta}=0.	
\end{equation*}
Since $p$ and $\theta$ are continuous across the discontinuity, the density will be as well. Hence, for $\varkappa\neq0$, these conditions can be rearranged into the following jump relations
\begin{equation*}
	\RH{\rho}=0, \quad \RH{u}=0\quad \RH{\eta}=0, \quad \RH{j}=0,
\end{equation*}
meaning that no contact discontinuities are admissible in this case. 
\begin{rmrk*}
	In the particular case where $\varkappa=0$, the system degenerates into Euler equations of compressible flows. In this case, it is rather clear that the evolution equation for $j$ is completely decoupled from the rest. Indeed, the variables $\rho$, $u$ and $E$ are completely determined from their equations alone. This means that, similarly to Euler equations, in the presence of a contact discontinuity, the temperature field will still be discontinuous across it. This may force $j$ to behave like a Dirac delta function locally in this point. This does not come out as surprising since for $\varkappa=0$, the total energy does not depend on $j$ and the system is only weakly hyperbolic even in one dimension.  
\end{rmrk*}

\subsection{Hugoniot curve for hyperbolic thermal conductivity}
Since $\calM\neq0$ is continuous across the shock, one can then rewrite equation \eqref{eq:RH_j} as
\begin{align*}
	\calM\RH{\frac{j}{\rho}}+\RH{\theta}=0,
\end{align*}
and the energy jump relation \eqref{eq:RH_E} as
\begin{equation}
	\RH{\calM\prn{\frac{\calM^2}{2\rho^2} + \varepsilon+\frac{p}{\rho}} + \frac{\varkappa^2}{2}\calM\left(\left(\frac{j}{\rho} + \frac{\theta}{\calM}\right)^2 -\frac{\theta^2}{\calM^2}\right) }=0.
	\label{eq:RH_E2}
\end{equation}
Then, it follows from equations (\ref{eq:RH_rho}, \ref{eq:RH_j}) that 
\begin{equation*}
	\RH{\frac{\varkappa^2}{2}\calM\left(\frac{j}{\rho} + \frac{\theta}{\calM}\right)^2} = 0,
\end{equation*}
so that the equation \eqref{eq:RH_E2} further simplifies into 
\begin{equation}
	\RH{\prn{\frac{\calM^2}{2\rho^2}+\frac{p}{\rho} + \varepsilon} - \frac{\varkappa^2}{2}\frac{\theta^2}{\calM^2} }=0.
	\label{eq:RH_E3}
\end{equation}
At this point, it is rather more practical to use the specific volume $\nu$ instead of the density $\rho$. It follows from equation \eqref{eq:RH_rhou} that $\calM^2$ can be expressed as
\begin{equation*}
	\calM^2 = -\frac{\RH{P}}{\RH{\nu}}.
\end{equation*}
This finally allows us to obtain  from equation \eqref{eq:RH_E3} the following consequence of the Rankine--Hugoniot conditions 	
\begin{equation*}
	\varepsilon_R-\varepsilon_L+\frac{1}{2}(p_R+p_L)(\nu_R-\nu_L)+\frac{\varkappa^2}{2}(\nu_R-\nu_L) \frac{\theta_R^2-\theta_L^2}{p_R-p_L}=0.
\end{equation*}
Let us formulate the following definition :
\begin{definition*}
	We call Hugoniot curve with center $(\nu_0,p_0)$ the curve in the $(\nu,p)$-plane defined as 
	\begin{equation}
		\varepsilon-\varepsilon_0+\frac{1}{2}(p+p_0)(\nu-\nu_0)+\frac{\varkappa^2}{2} \frac{\left(\theta^2-\theta_0^2\right)(\nu-\nu_0)}{(p-p_0)}=0.
		\label{eq:Hugoniot}
	\end{equation}
\end{definition*}
Through the shock, the jump relations \eqref{eq:RH} are supplemented by the Clausius-Duhem inequality
\begin{equation}
	\RH{\rho(u-\mathcal{D})\eta+\frac{\varkappa^2}{\rho}j} \geq 0.
	\label{eq:Clausius}
\end{equation}
If the subscript $0$ denotes the state ahead of the shock front, then this inequality can be rewritten as
\begin{equation}
	\eta-\eta_0 -\frac{\varkappa^2}{\calM^2}(\theta-\theta_0 )\geq0
	\label{eq:criterion_shock}
\end{equation} 
This inequality can be seen as a necessary condition of admissibility of shock waves for our system. It should be checked for each branch of the Hugoniot curve \eqref{eq:Hugoniot}. Note that for $\varkappa=0$ one recovers the usual entropy inequality for Euler equations of compressible fluids. 
\subsection{Properties of the Hugoniot curve}
In this proof, and for the sake of lightness, we shall use the \textit{prime} to denote total differentiation with respect to the specific volume $v$, that is for any variable $q$ we have 
\begin{equation*}
	q' = \frac{dq}{dv} 
\end{equation*} 
We denote the pressure function as a function of $v$ along the Hugoniot curve of center $(v_0,p_0)$ as
\begin{equation*}
	p= g(v,v_0,p_0).
\end{equation*}
Then, we abuse the notations to introduce 
\begin{equation*}
	\quad \theta(v) = \theta(g(v,v_0,p_0),v) ,\quad \eta(v)=\eta(g(v,v_0,p_0),v), \quad \Psi(v) = \eta(v) -\eta_0 -\frac{\varkappa^2}{\calM(v)^2}(\theta(v)-\theta_0 ).
\end{equation*}
In what follows, the arguments of these functions will be omitted. Remark that 
\begin{equation*}
	\Psi' = \eta' - \frac{\varkappa^2}{\calM^2} \theta' +\calG (\theta-\theta_0), \qquad {\rm where} \quad \calG= \frac{\varkappa^2}{\calM^4} \prn{\calM^2}' 
\end{equation*}
By taking the differential of  the Hugoniot curve equation \eqref{eq:Hugoniot} and making use of the Gibbs identity \eqref{eq:Gibbs}, we obtain 
\begin{equation*}
	\theta d\eta -\frac{1}{2}(p-p_0)dv +\frac{1}{2}(v-v_0)dp  - \frac{\varkappa^2}{\calM^2}\theta d\theta + \frac{1}{2}\calG \left(\theta^2-\theta_0^2\right)dv=0. 
\end{equation*}
Differentiation with respect to $v$ gives us
\begin{equation*}
	\theta\prn{ \eta'- \frac{\varkappa^2}{\calM^2}\theta' +  \prn{\theta-\theta_0}\calG} -  \frac{1}{2}\calG\prn{\theta-\theta_0}^2 -\frac{1}{2}(p-p_0) +\frac{1}{2}(v-v_0)g'  =0,  
\end{equation*}
and therefore
\begin{equation}
	\theta \Psi'-\frac{1}{2}(p-p_0) +\frac{1}{2}(v-v_0)g' -  \frac{1}{2}\calG\prn{\theta-\theta_0}^2 =0.
	\label{eq:der1}
\end{equation}
A second and third differentiations with respect to $v$ yield successively
\begin{equation}
	\theta' \Psi'+\theta \Psi'' +\frac{1}{2}(v-v_0)g'' -  \frac{1}{2}\calG'\prn{\theta-\theta_0}^2-  \calG\theta'\prn{\theta-\theta_0}  =0,
	\label{eq:der2}
\end{equation}
and 
\begin{align}
	&\theta'' \Psi'+ 2\theta' \Psi''+ \theta \Psi''' +\frac{1}{2}g''+\frac{1}{2}(v-v_0)g'''-  \frac{1}{2}\calG''\prn{\theta-\theta_0}^2-\prn{2\calG'\theta'+  \calG\theta''}\prn{\theta-\theta_0} -  \calG\theta'^2   =0.
	\label{eq:der3}
\end{align}
A series expansion of $p$ in the neighborhood of $v_0$ allows us  to write
\begin{equation*}
	\calG = -\varkappa^2\frac{v-v_0}{\prn{p-p_0}^2} \prn{g'-\frac{p-p_0}{v-v_0}}=\left.- \frac{\varkappa^2}{2}\frac{g''}{g'^2}\right\vert_{v=v_0} + {O}(v-v_0),  
\end{equation*}
and therefore, at the center of the Hugoniot curve, we have from equations \eqref{eq:der1} and \eqref{eq:der2} that 
\begin{equation*}
	\left.\Psi'\right\vert_{v=v_0}\!\!\!\!=0, \quad  \left.\Psi''\right\vert_{v=v_0}\!\!\!\!=0, 
\end{equation*}
and finally from equation \eqref{eq:der3} we obtain
\begin{equation*}
	\left.\Psi'''\right\vert_{v=v_0}\!\!\!\!   =\left.-\frac{1}{2\theta_0}\prn{1+\varkappa^2 \frac{\theta'^2}{g'^2}}g''\right\vert_{v=v_0}.
\end{equation*}
Notice that similarly to the Euler equations, $\Psi( v )$ has second order contact with the straight line $\Psi=0$.

\subsection{Analysis of an example: polytropic gas equation of state}
\label{subs:psi}
\begin{figure}[H]
	\center
	\includegraphics{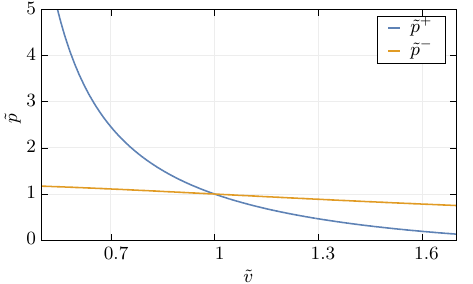}
	\caption{Plot of the functions $\tilde{p}^+$ and $\tilde{p}^-$ as a function of $\tilde v$.}
	\label{fig:p_tau}
\end{figure}

\begin{figure}[h]
	\includegraphics[width=0.5\textwidth]{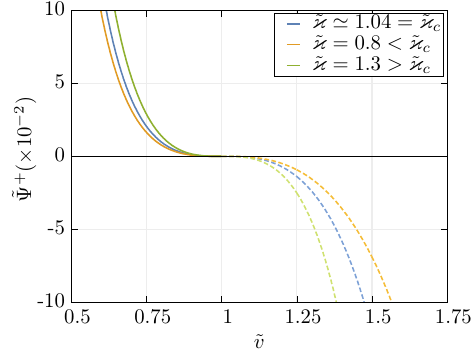} \includegraphics[width=0.5\textwidth]{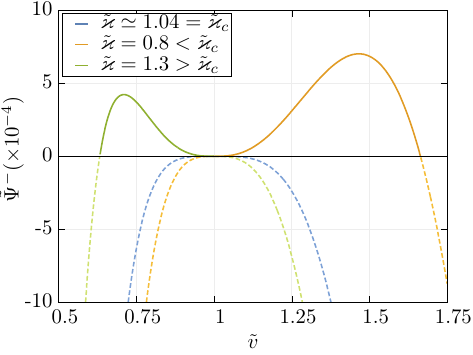}
	\caption{Plot of $\tilde\Psi^+$ (Left) and $\tilde\Psi^-$ (Right) as a function of $\tilde v $ for different values of $\tilde\varkappa$ that are $\varkappa_c= 1.0399...$, $\tilde\varkappa=0.8$ and $\tilde\varkappa=1.3$. $\gamma=2$. The dashed lines correspond to non-physical solutions that do not satisfy $\tilde\Psi\geq0$.}
	\label{fig:psi}
\end{figure}

\begin{figure}[H]
	\center
	\includegraphics[width=0.5\textwidth]{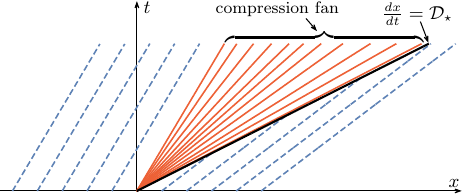}
	\caption{Schematic representation of the wave pattern in the $x-t$ plane, for the shock splitting solution. The shock propagates to the right, followed by a right facing compression fan. The Clausius-Duhem inequality is satisfied on the shock. The slope of the rightmost characteristic of the compression fan is equal to the velocity $\calD_\star$ of a shock  relating  the right  state  to the $\star$ state.}
	\label{fig:fan_xt}
\end{figure}

We carry on with the previous analysis for the specific example of a polytropic gas equation of state. In the scope of this study, casting the variables in dimensionless form allows to generalize the result so that it is independent from the choice of $\nu_0,p_0$. Therefore, let us consider the dimensionless specific volume $\tilde{\nu}$ and pressure $\tilde{p}$ defined as
\begin{equation*}
	\tilde{p} = p/p_0, \quad \tilde{v} = v/v_0,
\end{equation*}
and we introduce the dimensionless internal energy, temperature and entropy as follows
\begin{equation*}
	\tilde{\varepsilon} =  \frac{\tilde{v} \tilde{p}}{\gamma-1}, \quad	\tilde{\theta} =  \frac{\tilde{v} \tilde{p}}{\gamma-1}, \quad \tilde{\eta} = \mathrm{log}\prn{\tilde{p} \ \tilde{v}^\gamma}.
\end{equation*} 
One can verify that the Gibbs identity still holds in dimensionless variables. The Hugoniot curve can then be defined as
\begin{equation}
	\frac{\tilde\nu \tilde p-1}{\gamma-1}+\frac{1}{2}(\tilde p+1)(\tilde \nu-1)+\frac{\tilde\varkappa^2}{2} \frac{\left(\prn{\tilde \nu \tilde p}^2-1 \right)(\tilde \nu-1)}{(\tilde p-1)}=0, \quad \tilde\varkappa =  \frac{\nu_0 \ \varkappa}{(\gamma-1)c_V} 
	\label{eq:Hugoniot_nd}
\end{equation}
In order to further simplify the algebra and lighten the upcoming expressions, we will also take $\gamma=2$. Then, solving equation \eqref{eq:Hugoniot_nd} allows to obtain two distinct pressure functions $\tilde p^+(\tilde v)$ and $\tilde p^-(\tilde v)$, each corresponding to a different branch of the Hugoniot curve. The expressions of $\tilde p^+$ and $\tilde p^-$ can be obtained explicitly and are given by
\begin{equation*}
	\tilde p^{\pm}(\tilde v) = \frac{\tilde v+1\pm(1-\tilde v )\sqrt{\tilde\varkappa ^4 \tilde v ^2+\tilde\varkappa ^2 (\tilde v -1)^2+4}}{\tilde v  \left(\tilde\varkappa ^2 (\tilde v -1) \tilde v +3\right)-1}
\end{equation*}
Remark that in the limit $\tilde\varkappa=0$, $\tilde p^+$ degenerates towards the Euler equation's only admissible branch
\begin{equation*}
	\tilde p_E(\tilde v) = \frac{3-\tilde v}{3\tilde v-1}.
\end{equation*}
Thus, it seems reasonable here to call $\tilde{p}^+$ the \textit{acoustic} branch while we call $\tilde{p}^-$ the \textit{thermal} branch of the Hugoniot curve. A plot of these curves is given in figure \ref{fig:p_tau}. In comparison with Euler's equation, there exists a new branch $\tilde{p}^-$ due to the presence of thermal effects. 
Now we would like to check the sign of the quantity $\tilde\Psi$, as to examine the nature of admissible shock solutions for each of branches. Recall that those are shocks for which $\tilde\Psi\geq0$ along the Hugoniot curve. As shown in the previous paragraph, in the vicinity of the curve center, the positivity of $\tilde \Psi$ only depends on the sign of $\tilde{g}''(1)$ and whose expression as a function of $\tilde\varkappa$ is given hereafter by
\begin{equation*}
	\tilde g''(1)=\frac{2 \tilde\varkappa ^4}{\sqrt{\tilde\varkappa ^4+4}}-4 \tilde\varkappa ^2+\left(\tilde\varkappa ^2+3\right)^2-\left(\tilde\varkappa ^2+3\right) \left(\sqrt{\tilde\varkappa ^4+4}+1\right)
\end{equation*}
A brief study of this function in $\tilde\varkappa$ allows to see that it admits the non-trivial root $\tilde\varkappa_c= 1.0399...$ (in accordance with \eqref{eq:hyperplane}), such that $\tilde g''(1)<0$ for $\tilde\kappa<\tilde\kappa_c$ and $\tilde g''(1)>0$ for $\tilde\varkappa>\tilde\varkappa_c$. This implies that for $\tilde\kappa<\tilde\kappa_c$, there exists a neighborhood $B=]1,\tilde v_{max}[$ such that $\tilde\Psi>0\ \forall \; \tilde v \in B$. This leads to having a physically admissible branch of the Hugoniot curve in the region of lower density than the reference state \textit{i.e.}, the system can admit an expansion shock as a physically admissible solution. To further illustrate this, we provide in figure \ref{fig:psi} a graphical representation of $\tilde\Psi^+$ and $\tilde\Psi^-$, corresponding to the acoustic and thermal branches of the Hugoniot curve, respectively for different values of $\tilde\varkappa$.
The acoustic branch $\tilde\Psi^+$ shows a similar behavior to Euler's equation, and which is rather independent on $\tilde\varkappa$ in this regime. The main difference can be observed on the curves of $\tilde\Psi_-$, which shows a completely different behavior. Indeed, we can see here that depending on $\tilde\varkappa$, the system may admit additional shock wave solutions for $\tilde\varkappa>1$ or admit expansion shocks for $\tilde\varkappa<\tilde\varkappa_c$, in bounded regions around the center. Notice that for $\tilde\varkappa=\tilde\varkappa_c$ the system only admits shocks on the acoustic branch, meaning it behaves as Euler's equations. 
Finally, the last interesting property of our system that is related to the thermal branch of the Hugoniot curve is the so-called \textit{splitting phenomenon}. One can see on Figure \ref{fig:psi}  that $\Psi^-$ reaches its maximum at some point $\tilde v_\star$ both for $\tilde \varkappa <\tilde \varkappa_c$ and $\tilde \varkappa >\tilde \varkappa_c$. One can prove that at this point $\calM^2$ also attains its local  maximum \textit{i.e.},  $\prn{\calM^2}'=0,\;\prn{\calM^2}''< 0.$
Besides, for $\tilde v=\tilde v_\star$, $\left\vert u-\calD\right\vert$ is equal to the thermal sound speed $\sqrt{Y_1-Y_2}$. This is always the case of generic eigenfields that are neither genuinely nonlinear nor linearly degenerate  in the sense of Lax  \cite{liu1981admissible,lefloch2002hyperbolic}. If, for example,  $\tilde \varkappa >\tilde \varkappa_c$ and $\tilde  v _L<\tilde  v _\star$, one does not obtain for the corresponding Riemann problem the classical compression shock, but a rather compression shock followed by a compression fan (see Figure \ref{fig:fan_xt}). In this sense, the initial discontinuity verifying the corresponding Rankine-Hugoniot relations, splits into two distinct waves.
 The Clausius-Duhem inequality \eqref{eq:Clausius} is necessary but not sufficient to choose admissible shocks corresponding to the thermal branch. 
Indeed, admissible shocks must also verify the Oleinik-Liu admissibility criterion \cite{oleinik1959uniqueness,liu1981admissible}.  
This is also the case, for example, for the Euler equations of compressible fluids with non-convex adiabatic curves $p=g(v,\eta)$ (with constant $\eta$) \cite{wendroff1972riemann,wendroff1972riemann2,menikoff1989riemann,kluwick2018shock}. An example further illustrating this phenomenon will be shown later in the numerical tests.    
\section{Numerical scheme}
\label{sec:Numerics}

We provide in this section some numerical results for the model \eqref{eq:dissipative_hyp_heat_eq} which is written in one dimension of space in its conservative form as follows
\begin{equation}
	\pd{\mathbf{Q}}{t} + \pd{\mathbf{f}}{x} =\mathbf{S},
	\label{eq:PDE_num}
\end{equation}
where $\mathbf{Q}$, $\mathbf{f}$ and $\mathbf{S}$ are given by 
\begin{equation*}
	\mathbf{Q}=\prn{\begin{array}{c}
			\rho \\
			\rho u \\
			E \\
			j
	\end{array}}, \quad 	\mathbf{f}=\prn{\begin{array}{c}
			\rho u\\
			\rho u^2 + p \\
			E u + pu + \frac{\kappa^2}{\rho}\theta j \\
			ju + \theta
	\end{array}}, \quad 	\mathbf{S}=\prn{\begin{array}{c}
			0 \\
			0 \\
			0 \\
			-j/\tau
	\end{array}}.
\end{equation*}
In the case of a polytropic gas equation of state, $p$ and $\theta$ are computed as follows from the other variables
\begin{equation*}
	p = (\gamma-1)\prn{E-\frac{1}{2}\rho u^2-\frac{\kappa^2}{2\rho}j^2}, \quad \theta= \frac{p}{c_V (\gamma-1) \rho}
\end{equation*}
A numerical solution to this PDE will be sought on a discretization of the continuous domain $[x_L,x_R]$ comprised of $N$ uniform cells . The $i^{th}$ cell is delimited by $[x_{i-\frac{1}{2}},x_{i+\frac{1}{2}}]$. The numerical solution is evolved at discrete time-steps denoted by $t^n$. The cell centers $x_i$ and the discrete time $t^n$ are such that  
\begin{equation*}
	x_i=x_L+\left(i-\frac{1}{2}\right)\Delta x \quad \forall i\in[1,N] ,\quad \text{where} \ \Delta x = \frac{x_R-x_L}{N}, \qquad t^n=n\Delta t, \ n\in[0,n_{max}]
\end{equation*}
The hyperbolic system \eqref{eq:PDE_num} will be solved at the aid of a finite volumes method, based on a second-order ARS(2,2,2) scheme  \cite{Ascher1997}, which classifies as an Implicit-Explicit discretization (IMEX). We recall below the corresponding double butcher tableaux 
\begin{equation*}
	\begin{array}{c|ccc}
		0 & 0 & 0 & 0 \\
		\delta_1 & \delta_1 & 0 & 0 \\
		1 & \delta_2 & 1-\delta_2 & 0 \\
		\hline & \delta_2 & 1-\delta_2 & 0
	\end{array}, \quad \begin{array}{c|ccc}
		0 & 0 & 0 & 0 \\
		\delta_1 & 0 & \delta_1 & 0 \\
		1 & 0 & 1-\delta_1 & \delta_1 \\
		\hline & 0 & 1-\delta_1 & \delta_1
	\end{array}, \quad \left\{\begin{array}{l}
		\displaystyle \delta_1=1-\frac{1}{\sqrt{2}} \\[0.8em]
		\displaystyle  \delta_2=1-\frac{1}{2\delta_1}
	\end{array} \right. 
\end{equation*}
where the leftmost tableau gives the Runge-Kutta coefficients for the explicit part while the rightmost tableau gives those for the implicit part. In this context we choose to discretize the fluxes explicitly while we reserve the implicit discretization for the source terms, if any. As this is a two-stage scheme, at every time $t^n$, an intermediate stage $t^{n*}$ is first computed so that two-stage update formulas are given by 
\begin{alignat*}{2}
	&\mathbf{Q}_i^{n*} &&= \mathbf{Q}_i^n +  \Delta t\left( -\delta_1\frac{\mathbf{F}^n_{i+\frac{1}{2}}-\mathbf{F}^n_{i-\frac{1}{2}}}{\Delta x}+ \delta_1\mathbf{S}_i^{n*}\right) \\
	&\mathbf{Q}_i^{n+1} &&= \mathbf{Q}_i^n - \Delta t \left(\delta_2 \frac{\mathbf{F}^n_{i+\frac{1}{2}}-\mathbf{F}^n_{i-\frac{1}{2}}}{\Delta x} + (1-\delta_2)\frac{\mathbf{F}^{n*}_{i+\frac{1}{2}}-\mathbf{F}^{n*}_{i-\frac{1}{2}}}{\Delta x}  \right)+ \Delta t\left( (1-\delta_1) \mathbf{S}_i^{n*} + \delta_1 \mathbf{S}_i^{n+1}\right). 
\end{alignat*}   
In these expressions, $\mathbf{Q}_i^{n}$ are the cell averages and $\mathbf{S}_i^{n}$ are the corresponding source terms, evaluated in the $i^{th}$ cell at time $t^n$. $\mathbf{F}^n_{i+\frac{1}{2}}$ is the intercell flux approximated here by the Rusanov flux 
\begin{equation*}
	\mathbf{F}^n_{i+\frac{1}{2}} = \frac{1}{2}\prn{\mathbf{f}\prn{\mathbf{w_L}^{n}_{i+1}}+\mathbf{f}\prn{\mathbf{w_R}_{i}^{n}}} - \frac{1}{2}s_{i+\frac{1}{2}}^{max}\prn{\mathbf{w_L}^{n}_{i+1} - \mathbf{w_R}_{i}^{n}}
\end{equation*} 
where $\mathbf{w_L}_{i}^{n}$ (resp. $\mathbf{w_R}_{i}^{n}$) are the boundary extrapolated values to the left (resp. right) of $i^{th}$ cell, computed from a piece-wise linear reconstruction $\mathbf{w}_{i}^{n}$, constructed from $\mathbf{Q}_i^n$ by the MUSCL approach \cite{toro-book,VANLEER1977263}. $s_{i+\frac{1}{2}}^{max}$ is the maximum signal speed obtained as 
\begin{equation*}
	s_{i+\frac{1}{2}}^{max} = \max \prn{\max_{1\leq m \leq 4}\prn{\left\vert\lambda_m(\mathbf{w_L}_{i+1}^{n})\right\vert},\  \max_{1\leq m \leq 4}\prn{\left\vert\lambda_m(\mathbf{w_R}_{i}^{n})\right\vert} }
\end{equation*} 
As only the sources are implicit, a local fixed point iteration is employed to solve the resulting equation. The time-stepping of the scheme is subject to a CFL restriction such that at every time-step
\begin{equation*}
	\Delta t =\mathrm{CFL}  \frac{\Delta x}{\lambda_{max}}, \quad \lambda_{max} =\max_{i\in[1,N]} \prn{\max_{1\leq m \leq 4}\prn{\left\vert\lambda_m(\mathbf{w}_{i}^{n})\right\vert}}, 
\end{equation*}  
where the number $\mathrm{CFL}<1$ is a constant. Lastly, as numerical simulations of shock waves will be considered, a slope limiter will be used. In particular we will be using the Min-Mod slope limiter. 
\begin{rmrk*}
	When solving the original Euler system with heat conduction which is of mixed hyperbolic-parabolic type, the same scheme is used. The parabolic heat flux is simply incorporated to the explicit part by means of centered finite differences. The time-step in this case must also verify
	\begin{equation*}
		\Delta t < \frac{1}{2} \frac{\Delta x^2}{K}.
	\end{equation*} 
\end{rmrk*}
\section{Numerical results}
\label{sec:results}
In this part, whenever a test case is presented in dimensional variables, it will be assumed that every variable is in its corresponding units from the SI, \textit{i.e.}, $(Kg,m,s,K)$. While the units will be omitted from the text, they can still be seen on the plots for reference. In the description of initial data, be it for Euler's equation or our model, we shall report the used values of the density $\rho$, the velocity $u$ and the pressure $p$. In particular for our system, the value of $j$ will also be reported. Besides, in each case we will also give the value of the corresponding total energy, for the sake of completeness. In all of the cases presented below we use a polytropic gas equation of state.
\subsection{shock tube problem for Euler-equations with heat conduction}
\begin{figure}[H]
	\includegraphics[width=\textwidth]{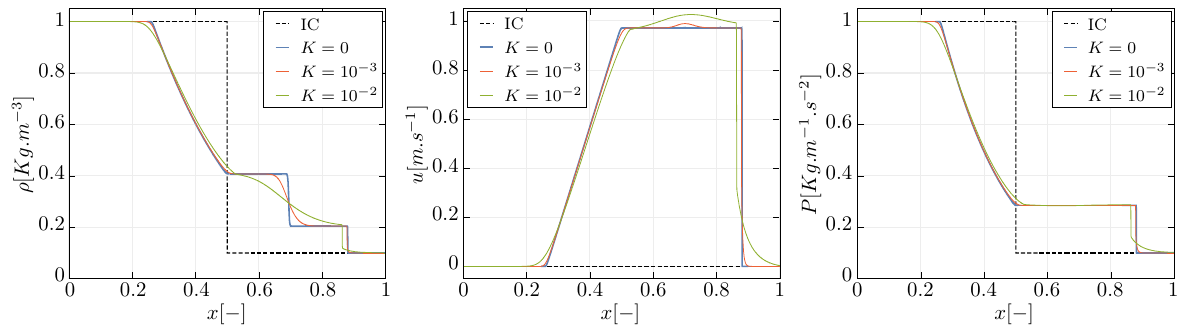}
	\caption{Numerical Profiles at of the density (Left), velocity (middle) and pressure (Right) obtained after solving numerically the Euler equations supplemented with heat conduction for the initial condition \eqref{eq:IC_choctube}. The solid line curves correspond to the numerical result at $t=0.2$ for different values of the heat conductivity $K=0$ (Blue), $K=10^{-3}$ (Red) and $K=10^{-2}$ (Green). The initial condition is plotted with black dashes lines. Here, $\gamma=1.4$ and $c_V =3/2$.}
	\label{fig:choc_tube}
\end{figure}

We consider here a classical shock tube problem. An initial discontinuity is placed at $x=x_s$, separating two different states on its left and right sides and which are given respectively by  
\begin{alignat}{3}
	&\rho_L = 1.0, \quad &&u_L = 0.0, \quad p_L = 1.0, \quad \prn{E_L = 2.50}   &&\qquad \text{for } x<x_s \nonumber \\
	&\rho_R = 0.1, \quad &&u_R = 0.0, \quad p_R = 0.1, \quad \prn{E_R = 0.25}  &&\qquad \text{for } x\geq x_s \label{eq:IC_choctube}
\end{alignat}   
In the context of the Euler equations, such an initial discontinuity gives rise to three distinct waves: a shock propagating to the right, a rarefaction propagating to the left and a contact discontinuity moving along the flow. In the case where heat conduction is also considered in the energy flux, the contact discontinuity is smeared out, leaving place to a smooth transition that is expanding over time. The numerical result is provided in figure \ref{fig:choc_tube}. The numerical domain is $[0,1]$ and the initial discontinuity position is $x_s=0.5$.
Note that for high heat conductivity, the overall profile of the solution begins to significantly change at later times due to the interaction of the diffused region with the shock front and the rarefaction fan. 
\subsection{Shock-tube problem for the hyperbolic model with relaxation}
In this part, we show numerical evidence of the asymptotic correspondence between the Euler equations supplemented with heat conduction with our model accounting for relaxation. Hence, let us consider the initial Riemann problem corresponding to the previous test case 
\begin{alignat}{3}
	&\rho_L = 1.0, \quad &&u_L = 0.0, \quad p_L = 1.0, \quad j_L=0.0, \quad \prn{E_L = 2.50}   &&\qquad \text{for } x<x_s \nonumber \\
	&\rho_R = 0.1, \quad &&u_R = 0.0, \quad p_R = 0.1, \quad j_R=0.0, \quad \prn{E_R = 0.25}  &&\qquad \text{for } x\geq x_s \label{eq:IC_choctube_Hyp}
\end{alignat} 

\begin{figure}[H]
	\includegraphics[width=\textwidth]{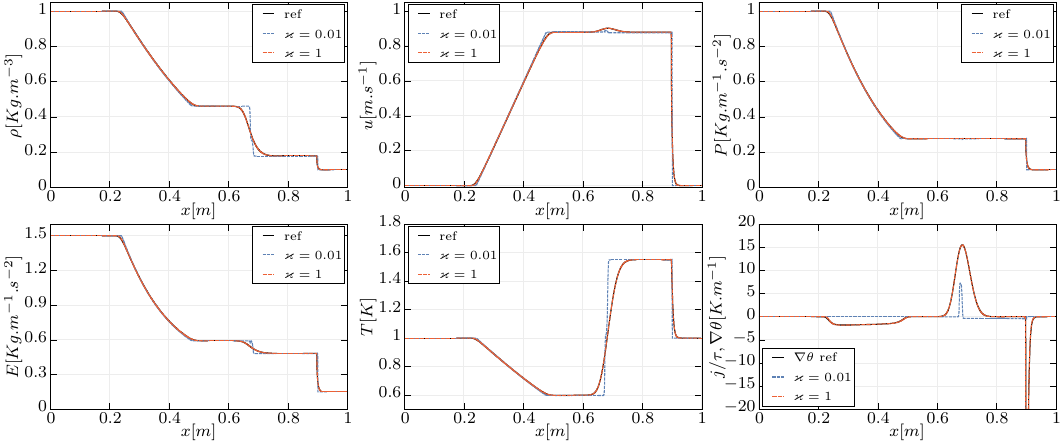}
	\caption{Numerical result for the choc tube problem \eqref{eq:IC_choctube_Hyp}. Simulation is obtained by the numerical scheme described in the previous section, on the computational domain $[0,1]$, discretized over $N=10000$ cells. The solid black line represents the reference solution obtained by solving the Euler system with heat conduction. The dashed lines correspond to the numerical solutions to the hyperbolic model with $\varkappa=0.01$ (Blue) and $\varkappa=1$ (Red). The solution  is given at final time $t=0.5$. $\mathrm{CFL}$ number is taken equal to $0.9$ (except for $\varkappa=1$ where $\mathrm{CFL}=0.5$). Values of the parameters are $\gamma=5/3$, $c_V=3/2$, $K=10^{-3}$. }
	\label{fig:choc_tube_hyp}
\end{figure}
For this test, we fix the value of the heat conductivity $K=10^{-3}$ and we set the relaxation time as defined from the asymptotic analysis relation \eqref{eq:tau} and whose expression in this case is recalled below
\begin{equation*}
	\tau = Kc_V(\gamma-1)\frac{\rho^2}{\varkappa^2p}.
\end{equation*}
We choose in this case to fix $\varkappa$ as a constant and have the relaxation time depend on the state variables as given by the above expression. In this setting, we show on figure \ref{fig:choc_tube_hyp} the numerical solution obtained for the initial data \eqref{eq:IC_choctube_Hyp} for two different values of $\varkappa$.
The numerical solution shows that for $\varkappa=1$, an excellent agreement of both models is observed. Note that the field $j/\tau$ also agrees with the temperature gradient, computed here with finite differences. For significantly low values of $\varkappa$, the behavior is rather closer to the classical Euler equations without any heat conduction.

\subsection{Expansion shock solution}
We have shown previously that expansion shock solutions are admissible for a value of $\varkappa$ above some critical value. Numerical evidence of such a solution will be presented in this part. The governing equations are taken here without relaxation terms and we consider an initial discontinuity for which we fix the right state   
\begin{equation}
	\rho_R = 1.0, \quad u_R = 0.0, \quad p_R = 1.0, \quad j_R=0.0, \quad \prn{E_L = 1.0}  \qquad \text{for } x \geq x_s.
	\label{eq:IC_RS1}
\end{equation}
According to our previous analysis, we need to have $\varkappa<\varkappa_c$ in order for expansion shocks to be admissible. We take $\gamma=2$ and $c_V=1$ and let us take the value $\tilde\varkappa=0.8$, for which the thermal branch was already plotted on figure \ref{fig:psi} and which satisfies this condition. From figure \ref{fig:psi}, we can pick for instance $\tilde  v =1.25$, to which corresponds the value $\rho_L=0.8$. Then the complete values for the left state are determined from the Rankine Hugoniot conditions and are as follows
\begin{equation}
	\rho_L = 0.8, \ u_L = -\frac{1}{4} \sqrt{1-\frac{\sqrt{13}}{5}}, \ p_L = \frac{3}{4}+\frac{\sqrt{13}}{20}, \ j_L=\frac{1}{4} \sqrt{\frac{1}{15} \left(11+\sqrt{13}\right)}, \ \prn{E_L = \frac{7}{150} \left(17+\sqrt{13}\right)}  \quad \text{for } x<x_s.
	\label{eq:IC_RS2}
\end{equation}
\begin{figure}[H]
	\includegraphics[width=\textwidth]{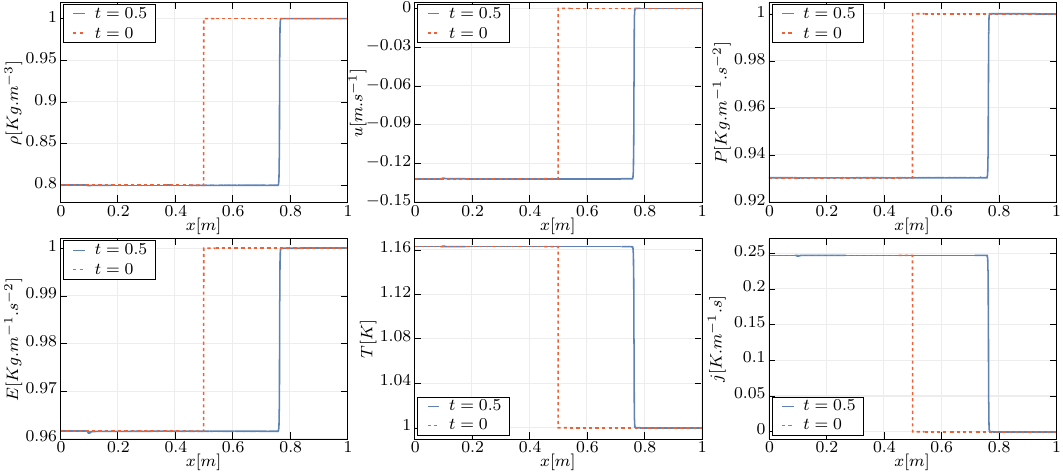}
	\caption{Numerical result for an expansion  shock solution with initial data (\ref{eq:IC_RS1}-\ref{eq:IC_RS2}). Simulation is obtained by the numerical scheme described in the previous section, on the computational domain $[0,1]$, discretized over $N=10000$ cells. The red dashed line corresponds to the initial condition while the solid blue line is the numerical solution at final time $t=0.5$. $\mathrm{CFL}$ number is taken equal to $0.9$. Values of the parameters are $\gamma=2$, $c_V=1$ and $\varkappa=0.8$. }
	\label{fig:rarefaction_shock}
\end{figure}
The numerical result on figure \ref{fig:rarefaction_shock} shows that the initial discontinuity propagates towards the region of higher density. Furthermore, notice that the profile of $u$ is such that the velocity is actually higher on the right side of the shock.

\subsection{Compression fan solution}
\begin{figure}[H]
	\includegraphics[width=\textwidth]{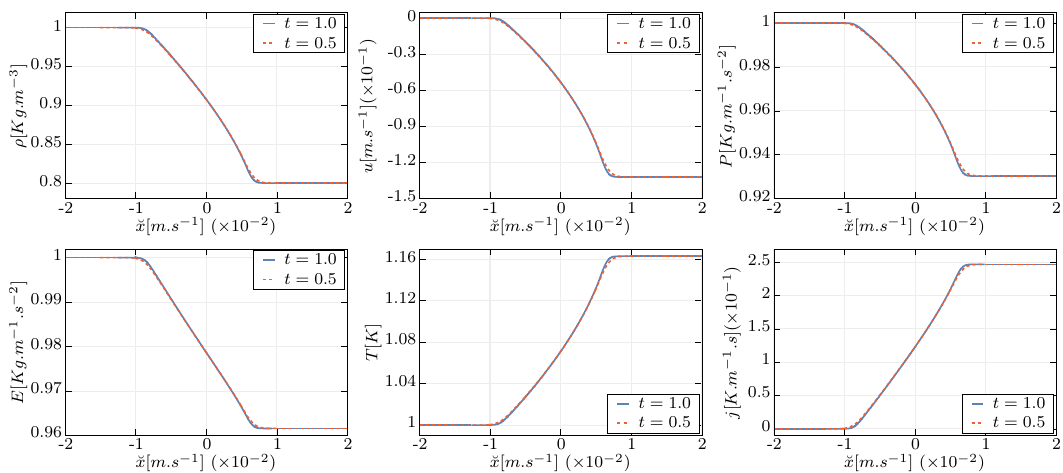}
	\caption{Numerical result for the compression fan solution with initial data (\ref{eq:IC_comp}), plotted as a function of $\breve{x}=(x-x_s-\calD t)/t$ on two different times $t=0.5$ (Red dashed line) and $t=1.0$ (Blue solid line). The full computational domain is $x\in[0,1]$, discretized over $N=100000$ cells. $\mathrm{CFL}$ number is taken equal to $0.9$. Values of the parameters are $\gamma=2$, $c_V=1$, $\varkappa=0.8$ and $x_s=0.2$. }
	\label{fig:simple_wave}
\end{figure}

Consider a self-similar solution to the governing equations \textit{i.e.}, a solution depending only on $x/t$. It is rather easy to formulate the ODE to find such a solution. However, it is even more straightforward to obtain it from the previous test case. To do this, one only needs to change the initial data of the previous problem by reversing the values on the left and right states so that they write 
\begin{gather}
	\rho_L = 1.0, \quad  u_L = 0.0, \quad p_L = 1.0, \quad j_L=0.0  \qquad \text{for } x<x_s. \nonumber \\
	\rho_R = 0.8, \ u_R = -\frac{1}{4} \sqrt{1-\frac{\sqrt{13}}{5}}, \ p_R = \frac{3}{4}+\frac{\sqrt{13}}{20}, \ j_R=\frac{1}{4} \sqrt{\frac{1}{15} \left(11+\sqrt{13}\right)}, \quad \text{for } x\geq x_s. 
	\label{eq:IC_comp}
\end{gather}
All the system parameters are naturally kept the same. It is well-known for the classical Euler equations, with a function $p(v,\eta)$, satisfying the inequalities $p_v<0, p_{vv}>0$, $p_\eta>0$ that such a reversal of initial data transforms  a shock into a rarefaction wave. It is very interesting to see that in our case, the previous expansion  shock solution turns into a compression fan. We take here a domain $[0,1]$ discretized over a refined mesh of $N=100000$ cells. The reason is that the solution, while smooth, spans over a narrow region of space. Thus, we would like to show that the regularity of the solution is not due to numerical effects but to its innate structure. The numerical solution is shown in figure \ref{fig:simple_wave}.
The solution is plotted as a function of the variable $\breve{x}=(x-x_s-\calD t)/t$ to emphasize the self-similar behavior. The wave speed $\calD$ is computed from the Rankine-Hugoniot relations and is given in this case by 
\begin{equation*}
	\calD = \sqrt{1-\sqrt{\frac{13}{5}}}.
\end{equation*}
A first glance at the figure might mislead into thinking this is a classical rarefaction wave but this is not the case. Indeed, notice from the profile of $u$ that the velocity to the right of the fan is lower than on the left of it.

\subsection{Shock wave splitting}
\begin{figure}[H]
	\includegraphics[width=\textwidth]{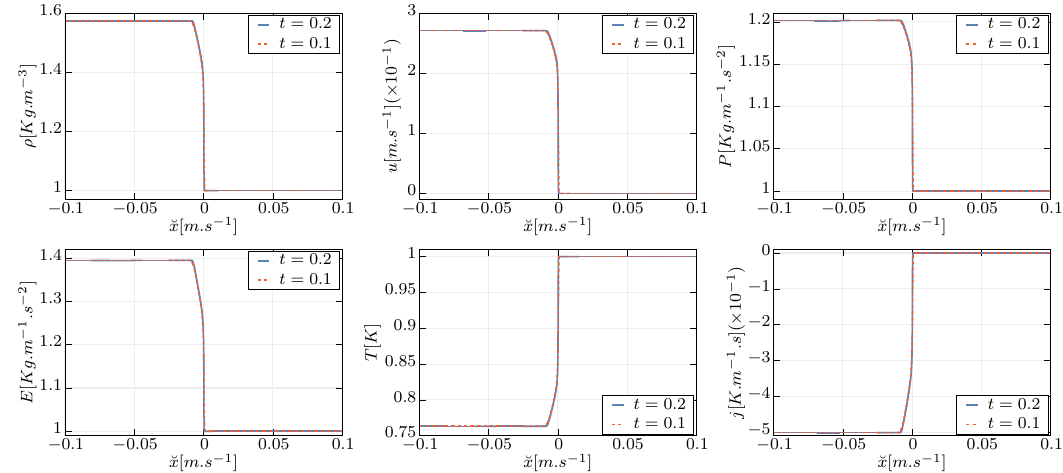}
	\caption{Numerical result for the shock wave splitting solution with initial data (\ref{eq:IC_split}), plotted as a function of the self-similar coordinate $\breve{x}=(x-x_s-\calD_\star t)/t$ on two different times $t=0.1$ (Red dashed line) and $t=0.2$ (Blue solid line). The full computational domain is $x\in[0,2]$, discretized over $N=100000$ cells. $\mathrm{CFL}$ number is taken equal to $0.5$. Values of the parameters are $\gamma=2$, $c_V=1$ and $\varkappa=1.3$. }
	\label{fig:splitting}
\end{figure}

For this test, we show numerical evidence of the shock splitting phenomenon, discussed in part \ref{subs:psi}. For this we consider $\varkappa=1.3>\varkappa_c$. We consider again the reference right state while the left state is computed from the Rankine-Hugoniot relation and Clausius-Duhem inequality. However, we emphasize that we take the left state such that the specific volume $v_L$ is on the left of the point $v_\star \simeq 0.7098$, where $\Psi_-$ reaches its maximum (see green curve of figure \ref{fig:psi}). We take for instance $v_L=0.635<v_\star$ , to which corresponds the following initial data
\begin{gather}
		\rho_L \simeq 1.575, \quad u_L \simeq 0.271, \quad p_L \simeq 1.202, \quad j_L\simeq-0.502, \quad \prn{E_L \simeq 1.395}  \qquad \text{for } x<x_s. \nonumber \\
	\rho_R = 1.0, \quad u_R = 0.0, \quad p_R = 1.0, \quad j_R=0.0, \quad \prn{E_L = 1.0}  \qquad \text{for } x\geq x_s.  
	\label{eq:IC_split}
\end{gather}
Approximate values of the left state variables were taken as to avoid their exact lengthy expressions. This shock is not admissible and it splits into a compression shock followed by a compression fan. The velocity of the shock $\calD_\star$ is computed in this case from the Rankine-Hugoniot conditions relating the right state and the $\star$ state, corresponding to the maximum of $\Psi_-$. The compression fan relates the star state with the left state. 
The numerical result is shown in figure \ref{fig:splitting}. We take here a domain $[0,2]$ with the a refined mesh as the previous test case. 

\section{Conclusion} 
\label{sec:conclusion}
We have presented in this paper a new first-order hyperbolic model describing heat transfer in a compressible fluid. The equations were derived from Hamilton's principle of stationary action. They can be written as a symmetric $t-$hyperbolic system in the sense of Friedrichs. The model was shown to be compatible with the classical Euler equations complemented by Fourier's heat conduction law when suitable relaxation terms are present in the governing equations. Evidence of expansion shocks as possible  solutions to the homogeneous system was established theoretically and confirmed via numerical simulations. This is due to the facts that the characteristic fields of thermal nature are neither genuinely non-linear nor linearly degenerate in the sense of Lax. Possible continuation of the analysis would be to prove rigorously the convergence of exact solutions of this model to those of the Euler system with heat conduction. At the numerical level, in order to solve multi-dimensional problems, it will be necessary to use proper numerical methods that are compatible with the curl-free constraint that binds the vector field $\mathbf{j}$. Such methods include for instance curl-cleaning \cite{MunzCleaning,Dedneretal,SHTCSurfaceTension,busto2021high,dhaouadi2022NSK,busto2023new} and exactly curl-free discretizations on staggered grids \cite{BalsaraSpicer1999,boscheri2021structure,SIGPR}.

%=============================================================================
%==========    A C K N O W L E D G M E N T S
\subsection*{Acknowledgments}
F.D. is a member of the INdAM GNCS group and acknowledges the financial support received from the Italian Ministry of Education, University and Research (MIUR) in the frame of the Departments of Excellence Initiative 2023--2027 attributed to DICAM of the University of Trento (grant L. 232/2016). F.D. was also funded by a \textit{UniTN starting grant} of the University of Trento and by the MUR grant of excellence \textit{PNRR Young Researchers, SOE}. The authors would also like to thank Michael Dumbser for important remarks and suggestions.

\newpage
\appendix
\section{Derivation of constrained Euler-Lagrange equations}
\label{app:calculus}
We first recall here the Lagrangian from which we shall recover the Euler-Lagrange equations 
\begin{equation*}
	\mathcal{L} = \int_{\Omega_t} \ L(\rho,\u,\gradphi,\dphi) \ d\Omega , \qquad L(\rho,\u,\gradphi,\dphi) = \frac{1}{2}\rho \nrm{\u}^2  - \frac{1}{2}\alpha(\rho)  \nrm{\gradphi}^2 - \rho \estar(\rho,\dphi)
\end{equation*}
In order to obtain the governing system of equations, one needs to consider two particular variations. A first variation yielding the least action with respect to the displacement of the continuum and a second variation for the independent degree of freedom $\phi$.
\subsection{Variation with respect to $\phi$}
The Euler-Lagrange equation corresponding to the variation of $\phi$  is 
\begin{equation*}
	-\pd{}{t}\prn{\pd{L}{\phi_t}} - \div{\pd{L}{\gradphi}} = 0,
\end{equation*}
which yields
\begin{equation}
	\pd{}{t}\prn{\rho \pd{\estar}{\dphi}} + \div{ \rho \pd{\estar}{\dphi}\u + \alpha(\rho) \gradphi} = 0
	\label{eq:EL_1}
\end{equation}

\subsection{Variation of the motion}
The Eulerian  variations (at fixed Eulerian coordinates) of the density and the velocity fields, under the total mass conservation constraint \eqref{eq:mass} can be obtained as a function of the virtual displacements  \cite{Berdichevsky2009,Gavrilyuk_2011,Dhaouadi2018,gavrilyuk2019dynamic}, and are given by 
\begin{equation*}
	{\delta} \rho = -\div{\rho \dx}, \quad {\delta} \u = \pd{\dx}{t} + \pd{\dx}{\x}\u - \pd{\u}{\x}\dx
\end{equation*}
Here $\dx(t, \mathbf x)$ is the virtual displacement of the  continuum. 
Hamilton's action between two arbitrary instants $t_1>t_0$ is 
\begin{equation*}
	a = \dblint{L(\rho,\u,\gradphi,\dphi)},
\end{equation*}  
thus yielding
\begin{equation*}
	\delta a = 0, \quad \text{or} \quad  \dblint{\delta L(\rho,\u,\gradphi,\dphi)}= \; 0.
\end{equation*}
Let us develop the latter expression further. 
\begin{align*}
	\dblint{\delta L} &= \dblint{\prn{\pd{L}{\rho}\delta \rho + \pd{L}{\u}\cdot\delta \u}} \\
	&=  \dblint{\prn{\frac{\nrm{\u}^2}{2} - \frac{1}{2}\alpha'(\rho) \nrm{\gradphi}^2 - \estar - \rho \pd{\estar}{\rho}}\delta \rho + \prn{\rho \u - \rho \pd{\estar}{\dphi}\gradphi}\cdot\delta \u}  \\
	& = \dblint{- \prn{\frac{\nrm{\u}^2}{2} - \frac{1}{2}\alpha'(\rho) \nrm{\gradphi}^2 - \estar - \rho \pd{\estar}{\rho}}\div{\rho \dx} 
		+  \prn{\rho \u - \rho \pd{\estar}{\dphi}\gradphi} \cdot \prn{\pd{\dx}{t} + \pd{\dx}{\x}\u - \pd{\u}{\x}\dx}}
\end{align*}
We recall, that by definition of the virtual displacement, $\dx$ is vanishing on the boundary $\partial\prn{\Omega\times[t_0,t_1]}$. Therefore, we can use Green-Ostrogradski's theorem to transform the previous equation as follows
\begin{align*}
	\dblint{\delta L}	& = \dblint{  \dx \cdot  \prn{\rho\nabla\prn{\frac{\nrm{\u}^2}{2}- \frac{1}{2}\alpha'(\rho) \nrm{\gradphi}^2 - \estar - \rho \pd{\estar}{\rho}}}}\\
	& + \dblint{ - \dx\cdot \pd{}{t} \prn{\rho \u - \rho \pd{\estar}{\dphi}\gradphi}+ \pd{\dx}{\x}\u \cdot \prn{\rho \u - \rho \pd{\estar}{\dphi}\gradphi}- \pd{\u}{\x}\dx \cdot\prn{\rho \u - \rho \pd{\estar}{\dphi}\gradphi} } 
\end{align*}
Straightforward calculus allows to further transform the products appearing in the latter integral as follows   
\begin{align*}
	\dblint{\pd{\dx}{\x}\u \cdot \prn{\rho \u - \rho \pd{\estar}{\dphi}\gradphi}} &= \dblint{\tr{\u \otimes \prn{\rho \u - \rho \pd{\estar}{\dphi}\gradphi} \ \pd{\dx}{\x} }} \\
	&= \dblint{- \div{\rho \u \otimes \u - \rho \u \otimes \pd{\estar}{\dphi}\gradphi} \cdot \dx} \\
\end{align*}
\begin{align*}
	&\dblint{\pd{\u}{\x}    \dx  \cdot  \prn{\rho \u - \rho \pd{\estar}{\dphi}\gradphi}}  = \dblint{\tr{\dx \otimes \prn{\rho \u - \rho \pd{\estar}{\dphi}\gradphi} \ \pd{\u}{\x}}}\\
	& = \dblint{-\div{\dx \otimes \prn{\rho \u - \rho \pd{\estar}{\dphi}\gradphi}} \cdot \u} 
	 = \dblint{-\div{\dx} \prn{\rho \u - \rho \pd{\estar}{\dphi}\gradphi} \cdot \u - \pd{}{\x}\prn{\rho \u - \rho \pd{\estar}{\dphi}\gradphi}\dx\cdot \u} \\
	& = \dblint{\grad{\prn{\rho \u - \rho \pd{\estar}{\dphi}\gradphi} \cdot \u}\cdot \dx - \grad{\rho \u - \rho \pd{\estar}{\dphi}\gradphi}\u \cdot \dx} = \dblint{\prn{\prn{\pd{\u}{\x}}^T\prn{\rho \u - \rho \pd{\estar}{\dphi}\gradphi}}  \cdot \dx}
\end{align*}
Replugging these expressions into the action allows to write it as follows
\begin{align*}
	\delta a = 	\int_{t_0}^{t_1} \!\!\! \int_\Omega \ \Bigg(\quad &\rho \nabla\prn{\frac{\nrm{\u}^2}{2}- \frac{1}{2}\alpha'(\rho) \nrm{\gradphi}^2 - \estar -  \rho \pd{\estar}{\rho}} - \pd{}{t} \prn{\rho \u - \rho \pd{\estar}{\dphi}\gradphi}  \\ 
	&-\div{\rho \u \otimes \u - \rho \u \otimes \pd{\estar}{\dphi}\gradphi}- \prn{\pd{\u}{\x}}^T\prn{\rho \u - \rho \pd{\estar}{\dphi}\gradphi}\quad \Bigg)\cdot \dx \ d\Omega \ dt = 0
\end{align*}
for any virtual displacement $\dx$, which allows us to recover the equation
\begin{equation*}
	\pd{}{t} \prn{\rho \u - \rho \pd{\estar}{\dphi}\gradphi} +\div{\rho \u \otimes \u - \rho \u \otimes \pd{\estar}{\dphi}\gradphi}-	\rho \nabla\prn{\frac{\nrm{\u}^2}{2}- \frac{1}{2}\alpha'(\rho) \nrm{\gradphi}^2 - \estar -  \rho \pd{\estar}{\rho}} + \prn{\pd{\u}{\x}}^T\prn{\rho \u - \rho \pd{\estar}{\dphi}\gradphi}=0.
\end{equation*}
At this point, we use the Euler-Lagrange equation \eqref{eq:EL_1} in order to recast the previous equation as
\begin{align*}
	&\pd{}{t} \prn{\rho \u} +\div{\rho \u \otimes \u +\alpha(\rho)\gradphi \otimes \gradphi}+	\rho \nabla\prn{-\frac{\nrm{\u}^2}{2}+ \frac{1}{2}\alpha'(\rho) \nrm{\gradphi}^2 + \estar +  \rho \pd{\estar}{\rho}} \\
	&+ \prn{\pd{\u}{\x}}^T\prn{\rho \u - \rho \pd{\estar}{\dphi}\gradphi}-\alpha(\rho) \pd{\gradphi}{\x} \gradphi - \rho \pd{\estar}{\dphi}\pd{\gradphi}{\x}  \u + \rho \pd{\estar}{\dphi} \grad{\gradphi\cdot \u -\dphi}=0
\end{align*}
Besides, we recall the following identities
\begin{align*}
	&\nabla\dot{\phi} =  \pd{\nabla \phi}{t} +  \pd{\gradphi}{\x} \u +  \prn{\pd{\u}{\x}}^T \gradphi \\
	&\grad{\frac{1}{2}(\rho\alpha'(\rho)-\alpha(\rho))\nrm{\gradphi}^2} = (\rho\alpha'(\rho)-\alpha(\rho))\pd{\gradphi}{\x} \gradphi +\frac{\nrm{\gradphi}^2}{2}\rho\alpha''(\rho)\nabla\rho = -\alpha(\rho)\pd{\gradphi}{\x} \gradphi + \rho \nabla \prn{ \frac{1}{2}\alpha'(\rho) \nrm{\gradphi}^2}  \\
	&   \rho\nabla\prn{  \estar + \rho \pd{\estar}{\rho}} = \rho \pd{\estar}{\dphi}\nabla \dphi + \rho \pd{\estar}{\rho}\nabla \rho + \rho\grad{\rho \pd{\estar}{\rho} } = \rho \pd{\estar}{\dphi}\nabla \dphi + \nabla \prn{\rho^2\pd{\estar}{\rho}}
\end{align*}
which finally leads us to the momentum equation in the form 
\begin{subequations}
	\begin{equation*}
		\frac{\partial \rho\u}{\partial t} + \div{ \rho \u \otimes \u + \prn{\rho^2\pd{\estar}{\rho}+\frac{1}{2}(\rho\alpha'(\rho)-\alpha(\rho))\nrm{\gradphi}^2}I + \alpha(\rho)\ \gradphi \otimes \gradphi }= 0.
	\end{equation*}
\end{subequations}

\section{Energy conservation and compatible curl-cleaning}
\label{app:energy_conservation}
In order to prove the energy conservation, it is generally more convenient to use conservative variables. We shall do the development in the case where arbitrary curl cleaning is also present in the equations. The total energy in this case is expressed by \eqref{eq:energy_cleaning}, written here for general $\alpha(\rho)$: 
\begin{equation*}
	\calE(\rho,\mathbf{m},\j,s,\ppsi) = \frac{1}{2\rho}\nrm{\mathbf{m}}^2 + \frac{1}{2}\alpha(\rho)\nrm{\j}^2 + \rho\varepsilon(\rho,s/\rho) + \frac{1}{2}\rho \nrm{\ppsi}^2,
\end{equation*}
thus implying that
\begin{equation*}
	\pd{\calE}{\mathbf{m}} = \frac{1}{\rho}\mathbf{m} = \u, \quad \pd{\calE}{\j} = \alpha(\rho)\j, \quad \pd{\calE}{s} = \theta, \quad P=\rho\pd{\calE}{\rho}+\mathbf{m}\cdot\pd{\calE}{\mathbf{m}}+s\pd{\calE}{s}-\calE 
\end{equation*}
Under these notations, we can rewrite our system of equations accounting for the relaxation sources and for curl-cleaning as follows 
\begin{alignat*}{2}
	& \frac{\partial \rho}{\partial t} &&+ \div{\rho \u } = 0 \\
	& \frac{\partial\mathbf{m}}{\partial t}  &&+ \div{ \mathbf{m} \otimes \u + \prn{\mathbf{m}\cdot\pd{\calE}{\mathbf{m}}+\rho\pd{\calE}{\rho}+s\pd{\calE}{s}-\calE} \ \mathbf{I}+\pd{\calE}{\j} \otimes \j }= 0 \\
	&\frac{\partial s}{\partial t} &&+ \div{s\u+ \pd{\calE}{\j}  } = \frac{\alpha (\rho)}{\theta(\rho,\eta)}\,\frac{\partial \calR}{\partial\j}\cdot  \j  \\
	& \frac{\partial \j}{\partial t} &&+ \frac{\partial \j}{\partial \x} \u + \prn{\frac{\partial \u}{\partial \x}}^T\j + \nabla\pd{\calE}{s}  + a_\j \nabla\wedge\ppsi= -\frac{\partial\calR}{\partial \j}  \\
	& \frac{\partial \ppsi}{\partial t} && + \pd{\ppsi}{\x}\u  - a_{\ppsi} \nabla \wedge \j= 0. 
\end{alignat*}
where arbitrary coefficients $a_\j$ and $a_{\ppsi}$ were introduced and whose exact expressions will be sought later. Multiplying each of the evolution equations by the corresponding conjugate variable, and summing up we obtain  
\begin{align*}
	& \pd{\calE}{\rho}\prn{\frac{\partial \rho}{\partial t}+\pd{\rho}{\x}\u + \rho\div{ \u }} +
	\pd{\calE}{\mathbf{m}}\cdot\prn{\frac{\partial\mathbf{m}}{\partial t} + \pd{\mathbf{m}}{\x}\u + \div{ \u}\mathbf{m}  + \nabla\prn{\mathbf{m}\cdot\pd{\calE}{\mathbf{m}}+\rho\pd{\calE}{\rho}+s\pd{\calE}{s}-\calE}+\div{\pd{\calE}{\j} \otimes \j }} \\
	+\ & \pd{\calE}{s}\prn{\frac{\partial s}{\partial t} +\pd{s}{\x}\u + s\div{\u}+ \div{\pd{\calE}{\j}  } } 
	+ \pd{\calE}{\j}\cdot\prn{\frac{\partial \j}{\partial t} + \frac{\partial \j}{\partial \x} \u + \prn{\frac{\partial \u}{\partial \x}}^T\j + \grad{\pd{\calE}{s}}  + a_\j \nabla\wedge\ppsi} \\
	+\ &\pd{\calE}{\ppsi}\cdot\prn{\frac{\partial \ppsi}{\partial t}  + \pd{\ppsi}{\x}\u  - a_{\ppsi} \nabla \wedge \j}= \pd{\calE}{s}\frac{\alpha (\rho)}{\theta(\rho,\eta)}\,\frac{\partial \calR}{\partial\j}\cdot  \j -\pd{\calE}{\j}\cdot\frac{\partial \calR}{\partial \j} = 0.
\end{align*}
Gathering the terms together yields
\begin{align*}
	& \pd{\calE}{t} + \pd{\calE}{\x}\u + \prn{\rho\pd{\calE}{\rho}+\mathbf{m}\cdot\pd{\calE}{\mathbf{m}}+s\pd{\calE}{s}} \div{ \u } +   \pd{}{\x}\prn{\rho\pd{\calE}{\rho}+\mathbf{m}\cdot\pd{\calE}{\mathbf{m}}+s\pd{\calE}{s}-\calE}\u+\div{\pd{\calE}{\j} \otimes \j }\cdot\u \\
	+\ & \pd{\calE}{s} \div{\pd{\calE}{\j}  }  + \grad{\pd{\calE}{s}}\cdot\pd{\calE}{\j}
	+ \pd{\calE}{\j}\cdot\prn{\frac{\partial \u}{\partial \x}^T\j   } -\pd{\calE}{\ppsi}\cdot\prn{ a_{\ppsi}\nabla \wedge \j}+ \pd{\calE}{\j}\cdot\prn{a_\j \nabla\wedge\ppsi}= 0
\end{align*}
Thus
\begin{align*}
	& \pd{\calE}{t} + \div{\prn{\rho\pd{\calE}{\rho}+\mathbf{m}\cdot\pd{\calE}{\mathbf{m}}+s\pd{\calE}{s}}\u +\pd{\calE}{s} \pd{\calE}{\j}} +\div{\pd{\calE}{\j} \otimes \j }\cdot\u
	+ \pd{\calE}{\j}\cdot\prn{\frac{\partial \u}{\partial \x}^T\j   } -\pd{\calE}{\ppsi}\cdot\prn{ a_{\ppsi} \nabla \wedge \j}+ \pd{\calE}{\j}\cdot\prn{a_\j \nabla\wedge\ppsi}= 0
\end{align*}
by identifying that 
\begin{equation*}
	\div{\prn{\pd{\calE}{\j} \otimes \j }\u} = 	\div{\pd{\calE}{\j} \otimes \j }\cdot\u + \tr{\prn{\pd{\calE}{\j} \otimes \j}\pd{\u}{\x}} = 	\div{\pd{\calE}{\j} \otimes \j }\cdot\u + \j\cdot\pd{\u}{\x} \pd{\calE}{\j} = \div{\pd{\calE}{\j} \otimes \j }\cdot\u + \pd{\calE}{\j}\cdot\prn{\frac{\partial \u}{\partial \x}^T\j   } 
\end{equation*}
one can finally write 
\begin{align*}
	& \pd{\calE}{t} + \div{\prn{\rho\pd{\calE}{\rho}+\mathbf{m}\cdot\pd{\calE}{\mathbf{m}}+s\pd{\calE}{s}}\u +\prn{\pd{\calE}{\j} \otimes \j }\u+\pd{\calE}{s} \pd{\calE}{\j}}  -a_{\ppsi}\pd{\calE}{\ppsi}\cdot\prn{  \nabla \wedge \j}+ a_\j\pd{\calE}{\j}\cdot\prn{ \nabla\wedge\ppsi}= 0
\end{align*}
The only terms remaining in non-conservative form are the ones linked to the curl cleaning. However, it is possible to cast it into a conservative form, provided a suitable choice of the coefficients $a_\j$ and $a_{\ppsi}$. In general, the following criteria is sufficient 
\begin{equation*}
	a_{\ppsi}\pd{\calE}{\ppsi}\cdot\j =  a_\j\pd{\calE}{\j}\cdot\ppsi = C \ \ppsi\cdot\j, \quad \forall\  \ppsi,\j \in \mathbb{R}^d 
\end{equation*} 
where $C$ is any real constant. In our case, replacing the energy derivatives by their expressions, this allows to write
\begin{equation*}
	a_{\ppsi} = \frac{C}{\rho}, \quad a_\j = \frac{C}{\alpha(\rho)}. 
\end{equation*} 
For the specific case $\alpha(\rho)=\varkappa^2/\rho$, we can take for instance $C=\varkappa a_c$ where $a_c$ is a free parameter then we have
\begin{equation*}
	a_{\ppsi} = \frac{\varkappa}{\rho}a_c, \quad a_\j = \frac{\rho}{\varkappa}a_c.
\end{equation*}
In this case, one obtains the energy conservation law 
\begin{equation*}
	\pdt{\calE} + \div{\ (\calE+P)\ \u +\prn{\alpha(\rho)\ \j \otimes \j}\ \u+\alpha(\rho)\, \theta(\rho,\eta)\,  \j + \varkappa a_c \j \wedge \ppsi} = 0.
\end{equation*}

\section{Symmetrization of the system}
\label{app:sym}
We recall that the variables $\mathbf{Q}$ and their conjugate $\mathbf{U}$ are given respectively by
\begin{equation*}
	\mathbf{Q}^T  = \prn{\rho,\m^T,s,\j^T} \quad \mathbf{U}^T = \prn{\mu, \u^T,\theta, \boldsymbol{\ell}^T}
\end{equation*}
In order to cast system of equations \eqref{eq:dissipative_hyp_heat_eq} into a symmetric $t-$hyperbolic from of Friedrichs, we first need to rearrange the terms in the evolution equation for $\j$ as follows
\begin{align*}
	\frac{\partial \j}{\partial t} +  \div{\mathbf{u}\otimes \j + \theta \  \mathbf{Id}} +  \prn{\frac{\partial \mathbf{u}}{\partial \x}}^T \j - \div{\mathbf{u}}\j = -\j/\tau
\end{align*} 
Then, it is required to alter the structure of the momentum equation by adding to it the term $-\prn{\nabla\wedge\j}\wedge E_\j \equiv 0$ so that it now writes 
\begin{equation*}
	\frac{\partial \m}{\partial t}  + \div{ \m \otimes E_\m + \j \otimes E_\j +\prn{p+\frac{1}{2}\prn{\rho\alpha'-\alpha}\j\cdot\j} \mathbf{I}} - \prn{ \frac{\partial \j}{\partial \x} - \prn{\frac{\partial \j}{\partial \x}}^T }E_\j  = 0 
\end{equation*}
Notice that the total pressure can be expressed as follows
\begin{align*}
	p  + \frac{1}{2}\prn{\rho\alpha'-\alpha} \nrm{\j}^2  = \rho \frac{\partial E}{\partial \rho}+\m\cdot E_\m+s E_H - E  = E^\star - \j \cdot E_{\j}, 
\end{align*}
thus allowing to further simplify the momentum equation as follows 
\begin{equation*}
	\frac{\partial \m}{\partial t}  + \div{ \m \otimes E_\m  +E^\star \mathbf{I}} +\div{E_\j}\j - \prn{\frac{\partial E_\j}{\partial \x}}^T \j= 0
\end{equation*}
These arrangements allow to write the following system of evolution equations 
	\begin{alignat*}{2}
		& \frac{\partial \rho}{\partial t} &&+ \div{\rho \mathbf{u} } = 0 \\
		& \frac{\partial \m}{\partial t}  &&+ \div{ \m \otimes E_\m  +E^\star \mathbf{I}} +\div{E_\j}\j - \prn{\frac{\partial E_\j}{\partial \x}}^T \j= 0 ,   \\
		& \frac{\partial \j}{\partial t} &&+  \div{\mathbf{u}\otimes \j + \theta \  \mathbf{Id}}- \div{\mathbf{u}}\j +  \prn{\frac{\partial \mathbf{u}}{\partial \x}}^T \j  = -\j/\tau.  \\
		&\frac{\partial s}{\partial t} &&+ \div{s\u+ E_\j  } = \frac{1}{\tau\theta}\,\j\cdot  E_\j 
	\end{alignat*}
We replace now $\mathbf{Q}^T = E^\star_\mathbf{U}$ and make use of the identities
\begin{equation*}
	\pd{\prn{E^\star \mathbf{u}}}{\mu} = E^\star_\mu \mathbf{u}, \quad    \pd{\prn{E^\star \mathbf{u}}}{\mathbf{u}} = \mathbf{u}\otimes E^\star_\mathbf{u} +E^\star   \mathbf{I}, \quad  \pd{\prn{E^\star \mathbf{u}}}{\boldsymbol{\ell}} = \mathbf{u}\otimes E^\star_{\boldsymbol{\ell}}, \quad   \pd{\prn{E^\star \mathbf{u}}}{\theta} = E^\star_\theta \mathbf{u}
\end{equation*}
in order to write the system as
	\begin{alignat*}{2}
		& \frac{\partial E^\star_\mu}{\partial t} &&+ \div{\pd{\prn{E^\star \mathbf{u}}}{\mu} } = 0 \\
		& \frac{\partial E^\star_\mathbf{u}}{\partial t}  &&+ \div{ \pd{\prn{E^\star \mathbf{u}}}{\mathbf{u}}}  + \div{\boldsymbol{\ell}}E^\star_{\boldsymbol{\ell}} - \prn{\pd{\boldsymbol{\ell}}{\x}}^T E^\star_{\boldsymbol{\ell}}= 0 , \quad  \\
		& \frac{\partial E^\star_{\boldsymbol{\ell}}}{\partial t} &&+  \div{\pd{\prn{E^\star \mathbf{u}}}{\boldsymbol{\ell}} }+\nabla\theta - \div{\mathbf{u}}E^\star_{\boldsymbol{\ell}}+  \prn{\frac{\partial \mathbf{u}}{\partial \x}}^T E^\star_{\boldsymbol{\ell}}  = -\j/\tau. \\
		&\frac{\partial E^\star_{\theta}}{\partial t} &&+ \div{\pd{\prn{E^\star \mathbf{u}}}{\theta} } +\div{ \boldsymbol{\ell} }= \frac{1}{\tau \theta}\,\j\cdot  E_\j 
	\end{alignat*}
which can be cast in compact notations as follows 
\begin{equation*}
	\tilde{\mathbf{A}}\pdt{\mathbf{U}} + \sum_i \prn{\tilde{\mathbf{B}}^i_1+\tilde{\mathbf{B}}^i_2} \pd{\mathbf{U}}{x_i}    = \mathbf{S}, 
\end{equation*}
where the matrices 
\begin{equation*}
	\tilde{\mathbf{A}} = \nabla^2_\mathbf{U}E^\star, \quad \tilde{\mathbf{B}}^i_1 = \nabla^2_\mathbf{U}(E^\star u_i)
\end{equation*}
are symmetric by definition, and $\tilde{\mathbf{B}}^i_2$ are the quasiliner matrices corresponding to the remaining non-conservative terms. A straightforward computation shows that they are also symmetric. For instance, in two space dimensions we have
\begin{equation*}
	\tilde{\mathbf{B}}^1_2 = \left(  
	\begin{array}{cccccc}
		0& 0 & 0 & 0 & 0 & 0  \\
		0& 0 & 0 & 0 & -E^\star_{\ell_2} & 0    \\
		0& 0 & 0 & E^\star_{\ell_2} & 0 & 0   \\
		0& 0 & E^\star_{\ell_2} & 0 & 0 & 1  \\
		0& -E^\star_{\ell_2} & 0 & 0 & 0 & 0    \\
		0& 0 & 0 & 1 & 0 & 0    
	\end{array}
	\right), \quad     \tilde{\mathbf{B}}^2_2 = \left(  
	\begin{array}{cccccc}
		0& 0 & 0 & 0 & 0 & 0  \\
		0& 0 & 0 & 0 & E^\star_{\ell_1} & 0    \\
		0& 0 & 0 & -E^\star_{\ell_1} & 0 & 0   \\
		0& 0 & -E^\star_{\ell_1} & 0 & 0 & 0  \\
		0& E^\star_{\ell_1} & 0 & 0 & 0 & 1    \\
		0& 0 & 0 & 0 & 1 & 0    
	\end{array}
	\right)
\end{equation*}
Therefore, if we consider the symmetric matrix $\tilde{\mathbf{B}}^i = \tilde{\mathbf{B}}^i_1+\tilde{\mathbf{B}}^i_2$, then system can be written in symmetric form 
\begin{equation*}
	\tilde{\mathbf{A}}\pdt{\mathbf{U}} + \sum_i \tilde{\mathbf{B}}^i \pd{\mathbf{U}}{x_i}    = \mathbf{S}. 
\end{equation*}

\bibliographystyle{mystyle}
\bibliography{./references}

\end{document}